\documentclass[12pt]{amsart}
\usepackage{amssymb}
\usepackage{amsmath}
\usepackage{amsthm}

\newcommand{\A}{{\mathcal A}}
\newcommand{\B}{{\mathcal B}}

\newcommand{\D}{{\mathcal D}}

\newcommand{\e}{\varepsilon}
\newcommand{\F}{{\mathcal F}}

\newcommand{\I}{{\mathcal I}}

\renewcommand{\P}{{\mathcal P}}

\newcommand{\p}{\partial}

\newcommand{\T}{{T^\ast M}}
\newcommand{\E}{{\mathcal E}}

\newtheorem{lemma}{Lemma}

\newtheorem{theorem}{Theorem}
\newtheorem{corollary}{Corollary}
\begin{document}
\title[Berezin-Toeplitz Quantization]
{A Formal Model of Berezin-Toeplitz Quantization}
\author[A.V. Karabegov]{Alexander V. Karabegov}
\address[Alexander V. Karabegov]{Department of Mathematics, Abilene Christian University, ACU Box 28012, Abilene, TX 79699-8012}
\email{axk02d@acu.edu}

\begin{abstract}
We give a new construction of symbols of the differential operators on the sections of a quantum line bundle $L$ over a K\"ahler manifold $M$ using the natural contravariant connection on $L$. These symbols are the functions on the tangent bundle $TM$ polynomial on fibres. 
 For high tensor powers of $L$, the asymptotics of the composition of these symbols leads to the star product of  a deformation quantization with separation of variables on $TM$ corresponding to some pseudo-K\"ahler structure on $TM$. Surprisingly, this star product is intimately related to the formal symplectic groupoid with separation of variables over $M$.
We extend the star product on $TM$ to generalized functions supported on the zero section of $TM$. The resulting algebra of generalized functions contains an idempotent element which can be thought of as a natural counterpart of the Bergman projection operator. Using this idempotent, we define an algebra of Toeplitz elements and show that it is naturally isomorphic to the algebra of Berezin-Toeplitz deformation quantization on $M$.
\end{abstract}

\date{June 10, 2006}
\maketitle

\section{Introduction}

Deformation quantization of a Poisson manifold $(M, \{\cdot,\cdot\})$ is an associative algebra structure on the space $C^\infty(M)[\nu^{-1},\nu]]$ of formal functions on $M$ with the product (named the star product)
\begin{equation}\label{E:star}
    \phi \star \psi = \sum_{r = 0} ^\infty \nu^r C_r(\phi,\psi),
\end{equation}
where $C_r$ are bidifferential operators on $M$ such that
\[
    C_0(\phi, \psi) = \phi\psi \mbox{ and } C_1(\phi,\psi) - C_1(\psi,\phi) = i\{\phi,\psi\}.
\]
Here we use the following definition of formal vectors. Given a vector space $V$, we denote by $V[\nu^{-1},\nu]]$ the space of formal Laurent series of the form
\[
    v = \sum_{r \geq n} \nu^r v_r,
\]
where $v_r \in V$ and $n$ is possibly negative, and call its elements formal vectors.
We assume that the unit constant 1 is the unity of the star algebra $(C^\infty(M)[\nu^{-1},\nu]], \star)$.
Two star products $\star_1$ and $\star_2$ are called equivalent if there exists a formal differential operator $B = 1 + \nu B_1 + \nu^2 B_2 + \ldots$ such that
\[
    \phi \star_1 \psi = B^{-1}(B\phi \star_2 B\psi)
\] 
for any formal functions $\phi,\psi\in C^\infty(M)[\nu^{-1},\nu]]$. A star product can be localized to any open subset $U \subset M$. 

Deformation quantization of Poisson manifolds was introduced in the seminal work \cite{BFFLS}. The existence and classification of deformation quantizations on the symplectic manifolds were established in a number of papers (see \cite{DWL}, \cite{F1},
\cite{OMY}, \cite{D}, \cite{F2}, \cite{NT}, 
\cite{BCG}, \cite{X}). Kontsevich showed the existence and gave the classification of deformation quantizations on an arbitrary Poisson manifold in \cite{K}. 

The concept of a star product is related to the notion of operator symbols and their composition. A star product $\phi \star \psi$ can be thought of as the asymptotic expansion of a family of symbol products $\phi *_h \psi$ depending on a small parameter $h$ as $h \to 0$, where the asymptotic parameter $h$ is replaced with the formal parameter $\nu$. For instance, this way one obtains Moyal deformation quantization from the composition of Weyl symbols. In some cases the space of $h$-dependent operator symbols does not carry a natural symbol product, but the symbol-operator correspondence has a well-defined classical limit as $h\to 0$ which allows to define a star product (see \cite{Guill},\cite{Sch}). The examples of star products related to symbols are mostly obtained from covariant and contravariant symbols on K\"ahler manifolds, introduced by Berezin (see \cite{Ber1}, \cite{Ber2}, \cite{CGR1}, \cite{CGR2}, \cite{Mor}, \cite{CMP2}, \cite{BMS}, \cite{Sch}, \cite{KSch}). These star products on K\"ahler manifolds enjoy a special property that the bidifferential operators $C_r$ in (\ref{E:star}) differentiate their first argument in holomorphic directions and the second argument in antiholomorphic ones or vice versa. The deformation quantizations with this property are called deformation quantizations with separation of variables or of the Wick type (see \cite{CMP1},\cite{BW},\cite{RT}). It was shown in \cite{CMP1} that all deformation quantizations with separation of
variables on a K\"ahler manifold can be explicitly constructed and bijectively parameterized by the formal deformations of the K\"ahler form (see also \cite{N}).

In this paper we consider a formal model of the construction of Berezin-Toeplitz quantization (see \cite{BMS}, \cite{Sch}). Berezin-Toeplitz quantization on a compact K\"ahler manifold $M$ uses the following data. Let $L$ be a quantum line bundle on $M$ (see the details in the main body of the paper). Denote by $\Pi^{(N)}$ the orthogonal projector  onto the space of holomorphic sections of 
$L^{\otimes N}$, the $N$-th tensor power of $L$ (the Bergman projector). To a given function $\phi \in C^\infty(M)$ there corresponds the Toeplitz operator $T^{(N)}_\phi = \Pi^{(N)} \circ \phi \circ \Pi^{(N)}$ on the sections of $L^{\otimes N}$ (here $\circ$ denotes composition of operators). The function $\phi$ is then called a contravariant symbol of the operator $T^{(N)}_\phi$. The symbol-operator mapping $\phi \mapsto T^{(N)}_\phi$ is not injective and therefore there is no natural product of contravariant symbols. However, one can extract a star product from the asymptotics of the Toeplitz operators as $N \to \infty$. It was proved in \cite{Sch} that there exists a unique deformation quantization (\ref{E:star}) on $M$, the Berezin-Toeplitz deformation quantization, such that for each $k$ 
\[
    ||T^{(N)}_\phi T^{(N)}_\psi - \sum_{r=0}^{k-1} \frac{1}{N^r}T^{(N)}_{C_r(\phi,\psi)}|| = O\left(\frac{1}{N^k}\right).
\]
In \cite{KSch} the Berezin-Toeplitz deformation quantization was completely identified. It was shown that the deformation quantization with the opposite star product is a deformation quantization with separation of variables whose characterizing formal deformation of the K\"ahler form was explicitly calculated.

In this paper we define symbols of the differential operators on the sections of $L^{\otimes N}$ and study the corresponding symbol product $*_{1/N}$. These symbols are the fibrewise polynomial functions on the tangent bundle $TM$. The asymptotic expansion of the symbol product $*_{1/N}$ as $N\to\infty$ leads to the star product $*$ of a deformation quantization with separation of variables on the tangent bundle $TM$ endowed with some pseudo-K\"ahler structure. It is not clear how to extend the symbols of differential operators introduced in this paper to wider classes of operators that would include, in particular, the Bergman projector $\Pi^{(N)}$. However, one can expect from semiclassical considerations that, as $N \to \infty$, such a symbol of $\Pi^{(N)}$ might have a well-defined limit which would be a generalized function supported at the zero section $Z$ of the tangent bundle $TM$. We show that the star product $*$ can be naturally extended to a class of generalized functions supported at $Z$ and that this class contains an idempotent element. Using this idempotent, one can define an algebra of Toeplitz elements. We show that this algebra is naturally isomorphic to the algebra of Berezin-Toeplitz deformation quantization on $M$.

\section{Deformation quantizations with separation of variables}

Let $M$ be a complex manifold endowed with a Poisson bivector field $\eta$ of type (1,1) with respect to the complex structure. We call such manifolds K\"ahler-Poisson. If $\eta$ is nondegenerate, $M$ is a pseudo-K\"ahler manifold. On a coordinate chart $U \subset M$ with local holomorphic coordinates $\{z^k, \bar z^l \}$, we write $\eta = ig^{\bar l k}\frac{\p}{\p z^k} 
\wedge \frac{\p}{\p\bar z^l}$. The condition that $\eta$ is Poisson is expressed in terms of the 
tensor $g^{\bar l k}$ as follows:
\begin{equation}\label{E:kp}
  g^{\bar l k} \frac{g^{\bar n m}}{\p z^k} = g^{\bar n k}  \frac{g^{\bar l m}}{\p z^k} \mbox{ and } g^{\bar l k} \frac{g^{\bar n m}}{\p\bar z^l} = g^{\bar l m} \frac{g^{\bar n k}}{\p\bar z^l}.
\end{equation}
The corresponding Poisson bracket  on $M$ is given locally by the formula
\begin{equation}\label{E:PoissM}
    \{\phi,\psi\}_M = ig^{\bar l k}\left(\frac{\p \phi}{\p z^k}\frac{\p\psi}{\p\bar z^l} - \frac{\p \psi}{\p z^k}\frac{\p\phi}{\p\bar z^l} \right).
\end{equation}
Deformation quantization (\ref{E:star}) on a K\"ahler-Poisson manifold $M$ is called a deformation quantization with separation of variables if the bidifferential operators $C_r$
differentiate their first argument in antiholomorphic directions and its second argument in holomorphic ones. 

Denote by $L_\phi$ the operator of star multiplication by a function $\phi$ from the left and by $R_\psi$ the operator of star multiplication by a function $\psi$ from the right, so that $L_\phi \psi = \phi \star \psi = R_\psi \phi$. Recall that the associativity of the star product $\star$ is equivalent to the condition that $[L_\phi,R_\psi]=0$ for any $\phi,\psi \in C^\infty(M)[\nu^{-1},\nu]]$.

With the assumption that the unit constant $1$ is the unity of the star product, the condition that $\star$ is a star product with separation of variables can be reformulated as follows. For any local holomorphic function $a$ and antiholomorphic function $b$ the operators $L_a$ and $R_b$ are the operators of pointwise multiplication by the functions $a$ and $b$, respectively, $L_a = a, \ R_b = b$. It is easy to check that 
\begin{equation}\label{E:C1}
  C_1(\phi,\psi) = g^{\bar l k} \frac{\p \phi}{\p \bar z^l}\frac{\p \psi}{\p z^k}.
\end{equation}
The Laplace-Beltrami operator $\Delta$ given locally by the formula 
\begin{equation}\label{E:Delta}
    \Delta = g^{\bar l k} \frac{\p^2}{\p z^k\p\bar z^l}
\end{equation}
is coordinate invariant and thus globally defined on $M$. 

For a given star product with separation of variables $\star$ on $M$ there exists a unique formal differential operator $B_\star$ on $M$ such that
\begin{equation}\label{E:ber}
B_\star(ab) = b\star a 
\end{equation}
for any local holomorphic function $a$ and antiholomorphic function $b$. The operator $B_\star$ is called the formal Berezin transform of the star product $\star$ (see \cite{Tr}). We see from Eqns. (\ref{E:C1}) and (\ref{E:Delta}) that 
\[
B_\star = 1 + \nu \Delta + \dots. 
\]
In particular, $B_\star$ is invertible. One can recover the star product of a deformation quantization with separation of variables from its formal Berezin transform. Using $B_\star$ as an equivalence operator, one can define a star product $\star'$ on $(M,\eta)$ as follows: 
\[
    \phi \star' \psi = B_\star^{-1}(B_\star\phi \star B_\star\psi).
\]
We call the star product opposite to $\star'$ the dual star product to the star product $\star$ and denote it by $\tilde \star$ so that
\[
    \phi \tilde\star \psi = \psi \star' \phi.
\]
It was shown in \cite{CMP3} that $\tilde \star$ defines a deformation quantization with separation of variables on the complex manifold $M$ endowed with the opposite Poisson bivector field $-\eta$, while the 
star product $\star'$ defines a deformation quantization with separation of variables on the manifold $\bar M$ with the opposite complex structure and the same Poisson bivector field $\eta$. The formal Berezin transform of the dual star product $\tilde\star$ is $B_\star^{-1}$ and the star product dual to $\tilde\star$ is $\star$. It follows from (\ref{E:ber}) that 
\begin{equation}\label{E:baabbb}
B_\star a = a \mbox{ and } B_\star b = b.
\end{equation}
In particular, $B_\star 1 = 1$. It was proved in \cite{CMP3} that for any local holomorphic function $a$ and antiholomorphic function $b$
\[
     B_\star aB_\star^{-1} = R_a \mbox{ and } B_\star bB_\star^{-1} = L_b.
\]

We call a formal differential operator $A = A_0 + \nu A_1 + \dots$ {\it natural} if the operator $A_r$ is of order not greater than $r$ for any $r \geq 0$. A star product (\ref{E:star}) is called natural if for each $r$ the bidifferential operator $C_r$ is of order not greater than $r$ with respect to each of its arguments (see \cite{GR}). For a natural star product the operators $L_f$ and $R_f$ are natural for any $f\in C^\infty(M)$. A star product with separation of variables on a K\"ahler-Poisson manifold is natural (see \cite{CMP3}, \cite{BW}, and \cite{N}).

Let $(M,\omega_{-1})$ be a pseudo-K\"ahler manifold. We say that a formal closed $(1,1)$-form $\omega = (1/\nu)\omega_{-1} + \omega_0 + \nu \omega_1 + \ldots$ is a formal deformation of the form 
$\omega_{-1}$. Deformation quantizations with separation of variables on $M$ are bijectively parametrized by the formal deformations of the pseudo-K\"ahler form $\omega_{-1}$ as follows.
Assume that $\omega$ is a formal deformation of the form $\omega_{-1}$. On a contractible coordinate chart $U \subset M$ the formal form $\omega$ has a potential $\Phi = (1/\nu)\Phi_{-1} + \Phi_0 + \nu \Phi_1 + \ldots$ such that $\omega = -i\p\bar \p \Phi$. There exists a unique star product on $(M, \omega_{-1})$ such that on each contractible coordinate chart $U$ for any holomorphic function $a$ and antiholomorphic function $b$ the following formulas hold:
\begin{eqnarray*}
   L_a = a,\ L_\frac{\p\Phi}{\p z^k} = \frac{\p}{\p z^k} + \frac{\p\Phi}{\p z^k} = e^{-\Phi} \left(\frac{\p}{\p z^k}\right) e^\Phi,\\
   R_b = b,\ R_\frac{\p\Phi}{\p\bar z^l} = \frac{\p}{\p\bar z^l} + \frac{\p\Phi}{\p\bar z^l} = e^{-\Phi} \left(\frac{\p}{\p \bar z^l}\right) e^\Phi.
\end{eqnarray*}
Namely, the centralizer of the operators $R_b$ and $R_{\p\Phi/\p\bar z^l}$ in the algebra of formal differential operators on $U$ can be identified with the algebra of left multiplication operators with respect to some star product $\star$ on $U$. This star product does not depend on the choice of the potential $\Phi$ and determines a global deformation quantization with separation of variables on $M$ parameterized by the formal form $\omega$.  

The star product $\tilde\star$ dual to $\star$ gives a deformation quantization with separation of variables on the pseudo-K\"ahler manifold $(M,-\omega_{-1})$. It is parameterized by a formal form $\tilde\omega = -(1/\nu)\omega_{-1} + \tilde\omega_0 + \nu \tilde\omega_1 + \ldots$.

A formal density $\rho = \sum_{r \geq n} \nu^r \rho_r$ (where $n$ is possibly negative) is called a trace density of a star product $\star$ on $M$ if for any functions $\phi,\psi \in C^\infty(M)$ such that at least one of them has compact support, the following identity holds:
\[
    \int_M \phi\star \psi\, \rho = \int_M \psi\star \phi\, \rho.
\]
Recall that a star product on a symplectic manifold has local derivations of the form $\delta = d/d\nu + A$, where $A$ is a formal differential operator (it does not contain derivatives with respect to the formal variable $\nu$). They are called $\nu$-derivations. On a symplectic manifold each star product has formal trace densities which differ by formal constant factors. There exists a canonical trace density $\mu$ uniquely determined by the folowing two requirements (see \cite{Tr}). 

\noindent (a) The leading term of the formal series $\mu$ is given by the formula
\begin{equation}\label{E:mu}
    \frac{1}{\nu^m m!}\omega_{-1}^m,
\end{equation}
where $m$ is the complex dimension of $M$. 

\noindent (b) Given any open subset $U\subset M$ with a $\nu$-derivation $\delta$ of the star product $\star$ on it and any function $f \in C^\infty_0 (U)$, the following identity holds:
\[
    \frac{d}{d\nu}\int_U f \, \mu = \int_U \delta(f) \, \mu.
\]
We call a formal function $\phi \in C^\infty(M)[\nu^{-1},\nu]]$ invertible if it can be represented as
\[
    \phi = \sum_{r \geq n} \nu^r\phi_r
\]
for some integer $n$ and with $\phi_n$ nonvanishing. Locally the function $\phi$ can be represented in the form
\[
     \phi = e^\theta,
\]
where
\[
    \theta = n \log \nu + \theta_0 + \nu \theta_1 + \ldots.
\]
Let $\star$ be the star product of the deformation quantization with separation of variables on a pseudo-K\"ahler manifold $(M,\omega_{-1})$ parameterized by a formal form $\omega$. 
In \cite{Tr} we gave an explicit construction of the canonical trace density of the star product $\star$. Denote by $B_\star$ the formal Berezin transform of the star product $\star$. Let $U$ be a contractible coordinate chart and $\Phi = (1/\nu)\Phi_{-1} +\ldots$ a formal potential of $\omega$ on $U$. There exists a potential $\Psi$ of the dual form $\tilde\omega$ expressed as
\[
      \Psi= -m \log \nu - (1/\nu)\Phi_{-1} + \Psi_0 + \nu \Psi_1 +  \ldots
\] 
and such that
\begin{equation}\label{E:normalization}
   B_\star\left(\frac{\p \Psi}{\p z^k}\right) = -\frac{\p \Phi}{\p z^k},\ B_\star\left(\frac{\p \Psi}{\p \bar z^l}\right) = -\frac{\p \Phi}{\p\bar z^l} \mbox{ and } B_\star\left(\frac{d \Psi}{d\nu}\right) = -\frac{d\Phi}{d\nu}. 
\end{equation}
The potential $\Psi$ is determined by Eqn. (\ref{E:normalization}) uniquely up to a constant summand. The first two equations in (\ref{E:normalization}) determine $\Psi$ only up to a formal constant. The canonical trace density of the star product $\star$ can be expressed on $U$ as follows:
\begin{equation}\label{E:tracedensity}
      \mu = C e^{(\Phi + \Psi)}dzd\bar z,
\end{equation}
where the constant factor $C$ is uniquely determined by requirement (\ref{E:mu}). 
The star products $\star, \star'$, and $\tilde\star$ have the same trace densities and the same canonical trace density. 

\section{Symbols of differential operators on a quantum line bundle}\label{S:symb}

In this section we will relate to a pseudo-K\"ahler manifold $M$ endowed with a pseudo-K\"ahler form $\omega_{-1}$ a family of associative products on the fibrewise polynomial functions on the tangent bundle $TM$. Let $U$ be a contractible coordinate chart on $M$ with holomorphic coordinates $\{z^k,\bar z^l\}$. Consider a holomorphic Hermitian line bundle $L$ on $U$ with the fibre metric $|\cdot|$ whose metric-preserving covariant connection has the curvature $\omega_{-1}$. Such a line bundle is called `quantum'. Let $s$ be a nonvanishing holomorphic section of 
$L$. Then $\Phi_{-1} = - \log |s|^2$ is a potential of the form $\omega_{-1}$ so that $\omega_{-1} = -i\p\bar \p \Phi_{-1}$. In coordinates, $\omega_{-1} = -ig_{k\bar l} dz^k \wedge d\bar z^l$, where $g_{k\bar l} = \frac{\p^2 \Phi_{-1}}{\p z^k \p \bar z^l}$. The inverse matrix $(g^{\bar lk})$ of $(g_{k\bar l})$ determines a (global) Poisson bivector field $\eta = ig^{\bar lk} \frac{\p}{\p z^k} \wedge\frac{\p}{\p \bar z^l}$ and satisfies Eqn. (\ref{E:kp}).
For a given positive integer $N$ we will set
\[
                      h = \frac{1}{N}.
\]
Consider the $N$-th tensor power $L^{\otimes N}$ of the line bundle $L$ and denote
\[
          \Psi_{h} = \frac{1}{h}\, \Phi_{-1}.
\]
In the trivialization of $L^{\otimes N}$ determined by the section $s^N$
the metric-preserving covariant connection $\nabla_\bullet$ is as follows:
\begin{equation}\label{E:covar}
    \nabla_k = \frac{\p}{\p z^k} - \frac{\p \Psi_{h}}{\p z^k},\quad \nabla_{\bar l} = \frac{\p}{\p \bar z^l}.
\end{equation}
Introduce a contravariant connection $\nabla^\bullet$ on $L^{\otimes N}$ by lifting the index in (\ref{E:covar}) via the tensor $h g^{\bar lk}$:
\[
   \nabla^p = - h g^{\bar lp}\frac{\p}{\p \bar z^l}, \quad \nabla^{\bar q} = h g^{\bar qk} \left(\frac{\p}{\p z^k} - \frac{\p \Psi_{h}}{\p z^k}\right).
\]
Let $\{\eta^k,\bar\eta^l\}$ be the fibre coordinates on the tangent bundle $TU$ corresponding to the base coordinates $\{z^k,\bar z^l\}$. For a function $f=f(z,\bar z)$ on $U$ denote the pointwise multiplication operator by $f$ by the same symbol. Introduce the following operators on the sections of $L^{\otimes N}$:
\[
   \hat f = f,\ \widehat{\eta^p} = \nabla^p,\ \widehat{ \bar\eta^q} = \nabla^{\bar q}.
\]
We make a crucial observation that, according to Eqn. (\ref{E:kp}), the operators $z^1,\ldots,z^m, \nabla^1, \ldots, \nabla^m$ pairwise commute and the operators ${\bar z}^1,\ldots,{\bar z}^m, \nabla^{\bar 1}, \ldots, \nabla^{\bar m}$ pairwise commute as well. We can treat these commuting families as `coordinate' and `momentum' operators on the sections of $L^{\otimes N}$ and define symbols of differential operators on $L^{\otimes N}$ via the normal ordering. 
Denote by $\P_r(TU)$ the space of fibrewise polynomial functions on the tangent bundle $TU$ whose fibrewise degree is not greater than $r$ and set $\P(TU) = \cup_r \P_r(TU)$.
To a fibrewise polynomial function
\[
    P = \sum_i u_i(\eta)\, f_i\,  v_i(\bar\eta)
\]
from $\P(TU)$, where $u_i$ and $v_i$ are monomials in the fibre variables $\{\eta^k\}$ and $\{\bar\eta^l\}$, respectively, we relate a differential operator $\hat P$ on $L^{\otimes N}$ represented in the normal form
\[
   \hat P = \sum_i u_i(\hat{\eta})\, f_i\,  v_i(\widehat{\bar\eta}).
\]

The symbol-operator mappings $P\mapsto \hat P$ is a bijection of the space of fibrewise polynomial functions on $TU$ onto the space of differential operators on the sections of the line bundle $L^{\otimes N}$. The symbol product will be denoted $*_{h}$ so that for symbols $P,Q \in \P(TU)$
\[
   \widehat{(P *_{h} Q)} = \hat P \hat Q.
\]

Denote by $L_P$ and $R_P$ the operators of left and right multiplication by a symbol $P$ with respect to the product $*_{h}$, respectively. Notice that for functions $\phi,\psi$ on $U$ and monomials $u = u(\eta),v = v(\bar\eta)$ 
\begin{equation}\label{E:separ}
\phi *_{h} \psi = \phi\psi, \ L_u = u(\eta), \mbox{ and } R_v = v(\bar\eta).
\end{equation}
Due to the commutation relations
\begin{equation}\label{E:commut}
[\widehat{\eta^p},\hat f] = - h\, g^{\bar lp}\frac{\p f}{\p \bar z^l} \mbox{ and }
 [\widehat{\bar\eta^q},\hat f] = h\, g^{\bar qk}\frac{\p f}{\p z^k}, 
\end{equation}
where $f = f(z,\bar z)$ is a function on $U$, we get that 
\begin{equation}\label{E:etaf}
f *_{h} \eta^p = \eta^p f + h\, g^{\bar lp}\frac{\p f}{\p\bar z^l} \mbox{ and }
\bar\eta^q *_{h} f = f\bar\eta^q + h\, g^{\bar qk}\frac{\p f}{\p z^k}.
\end{equation}
It follows from (\ref{E:separ}) and (\ref{E:etaf}) that the symbols which do not depend on the antiholomorphic fibre variables $\{\bar\eta^l\}$ form a subalgebra of the algebra $(\P(TU),*_{h})$ which will be denoted $\A$. Denote by $L_P^{\A}, R_P^{\A}$ the operators of left and right multiplication by a symbol $P$ in the algebra $\A$, respectively. Similarly, $\B$ will denote the subalgebra of symbols which do not depend on the variables $\{\eta^k\}$ and $L_P^{\B}, R_P^{\B}$ the left and right multiplication operators by a symbol $P$ in the algebra $\B$. For a function $f$ on $U$ 
\begin{equation}\label{E:raflbf}
       R^{\A}_f = f \mbox{ and } L^{\B}_f = f.   
\end{equation}
Formulas (\ref{E:separ}) and (\ref{E:etaf}) imply that 
\begin{equation}\label{E:ralb}
   R^{\A}_{\eta^p} =  \eta^p + h\, g^{\bar lp}\frac{\p }{\p\bar z^l}\mbox{ and } 
   L^{\B}_{\bar\eta^q} = \bar\eta^q + h\, g^{\bar qk}\frac{\p }{\p z^k}.
\end{equation}
Let $\D(TU)$ denote the space of differential operators on $TU$. It has a decreasing filtration by the subspaces $\D_r(TU)\subset \D(TU)$ of differential operators annihilating $\P_r(TU)$. Denote by 
$\hat \D(TU)$ the completion of the space $\D(TU)$ with respect to this filtration. The elements of $\hat \D(TU)$ act on the space $\P(TU)$ as differential operators of infinite order. Introduce operators $J_{h},K_{h}\in \hat\D(TU)$ by the formulas
\[
    J_{h} = \exp\left(h g^{\bar lk}\frac{\p^2}{\p \eta^k \p\bar z^l}\right) \mbox{ and } K_{h} = \exp\left(h g^{\bar lk}\frac{\p^2}{\p z^k \p\bar \eta^l}\right).
\]
Here the exponentials are given by convergent series in the topology determined by the filtration.

We need the following lemma. 
\begin{lemma}\label{L:commute}
 For any functions $\phi,\psi$ on $U$ the operators $J_{h}\phi J_{h}^{-1}$ and  $K_{h}\phi K_{h}^{-1}$ commute with the operator of point-wise multiplication by the function $\psi$,
\[
    \left[J_{h}\phi J_{h}^{-1}, \psi\right] = \left[K_{h}\phi K_{h}^{-1}, \psi\right] = 0. 
\]
\end{lemma} 
{\it Proof.} Denote by ${\mathcal S}\subset \hat\D(TU)$ the ring of operators  of the form $A = A(z,\bar z, \frac{\p}{\p \eta})$, i.e., generated by the multiplication operators by the functions $f = f(z,\bar z)$ and the derivations $\frac{\p}{\p\eta^k}$. The inclusion
\[
          \left[h g^{\bar lk}\frac{\p^2}{\p\eta^k \p\bar z^l}, {\mathcal S}\right] \subset {\mathcal S}
\]
implies that $J_{h}{\mathcal S}J_{h}^{-1}\subset {\mathcal S}$. The statement that the operators $J_{h}\phi J_{h}^{-1}$ and $\psi$ commute follows from the fact that the ring ${\mathcal S}$ is commutative and contains the multiplication operators by the functions $f = f(z,\bar z)$. The rest of the Lemma can be checked similarly.

Notice that for a function $f=f(z,\bar z)$  and monomials $u=u(\eta),v=v(\bar\eta)$
\begin{eqnarray}\label{E:thru}
   J_{h}f = J_{h}^{-1}f = K_{h}f = K_{h}^{-1}f = f,\\
 J_{h}u = J_{h}^{-1}u = u,\  K_{h}v = K_{h}^{-1}v = v.\nonumber
\end{eqnarray}
Using Eqns.(\ref{E:kp}) and (\ref{E:ralb}) one can check that
\begin{equation}\label{E:jetaj}
     R^{\A}_{\eta^p} = J_{h} \eta^p J_{h}^{-1}  \mbox{ and } L^{\B}_{\bar\eta^q} = K_{h} \bar\eta^q K_{h}^{-1}.  
\end{equation}
Taking into account Eqns. ({\ref{E:thru}}) and (\ref{E:jetaj}) we get
\begin{eqnarray}\label{E:important}
    f *_{h} u(\eta) = R^\A_{u}f = u\left(R^\A_{\eta}\right)f =\\ 
J_{h} u(\eta) J_{h}^{-1} f = J_{h}(uf) = J_{h}fJ_{h}^{-1} u. \nonumber
\end{eqnarray}
This calculation allows to determine the operator $L_f$.
\begin{lemma}\label{L:lfrf}
The operators of left and right multiplication by a function $f = f(z,\bar z)$ in the algebra $(\P(TU),*_{h})$ can be expressed as follows:
\[
    L_f = J_{h}fJ_{h}^{-1} \mbox{ and } R_f = K_{h}fK_{h}^{-1}.
\]
\end{lemma}
{\it Proof.} We will prove the first formula using Eqns. (\ref{E:separ}), (\ref{E:raflbf}), and (\ref{E:important}). For a function $\phi=\phi(z,\bar z)$ and monomials $u=u(\eta),v=v(\bar\eta)$ we get, taking into account Lemma \ref{L:commute} and that the operator $J_{h}$ commutes with the operator of pointwise multiplication by $v(\bar\eta)$, that
\begin{eqnarray*}
   f *_{h} (u\phi v) = ((f *_{h} u)*_{h}\phi)*_{h}v =\\
  \left(J_{h}fJ_{h}^{-1}u\right)\phi v = J_{h}fJ_{h}^{-1}(u\phi v).
\end{eqnarray*}
The second formula in the Lemma can be checked similarly.

Introduce a potential $\Xi_{h}$ on $U$ by the following formula:
\begin{equation}\label{E:potxi}
   \Xi_{h} = \Psi_{h} + \frac{\p\Psi_{h}}{\p z^k}\eta^k + \frac{\p\Psi_{h}}{\p\bar z^l}\bar\eta^l + \log g,
\end{equation}
where $g = \det(g_{k\bar l})$. We need to calculate a number of  commutation relations in the algebra of differential operators on $L^{\otimes N}$.
\begin{lemma}\label{L:commrel}
 The following formulas and commutation relations hold:
\begin{equation}\label{E:hatxi}
    \widehat{\frac{\p\Xi_{h}}{\p \eta^p}} = \frac{\p\Psi_{h}}{\p z^p},\ \widehat{\frac{\p\Xi_{h}}{\p z^p}} = \frac{\p}{\p z^p} - h g^{\bar lk} \frac{\p^2\Psi_{h}}{\p z^k \p z^p}\frac{\p}{\p\bar z^l};
\end{equation}
\begin{eqnarray}\label{E:hatxibar}
\widehat{\frac{\p\Xi_{h}}{\p \bar\eta^q}} = \frac{\p\Psi_{h}}{\p \bar z^q},\ \widehat{\frac{\p\Xi_{h}}{\p \bar z^q}} = - \frac{\p}{\p \bar z^q} +  h g^{\bar lk} \frac{\p^2\Psi_{h}}{\p \bar z^l \p \bar z^q}\frac{\p}{\p z^k} + \\
\frac{\p\Psi_{h}}{\p \bar z^q} - h g^{\bar lk} \frac{\p^2\Psi_{h}}{\p \bar z^l \p \bar z^q}\frac{\p\Psi_{h}}{\p z^k};\nonumber
\end{eqnarray}
\[
\left[\widehat{\frac{\p\Xi_{h}}{\p \eta^p}}, \widehat{\eta^k}\right] = \delta^k_p,\  \left[\widehat{\frac{\p\Xi_{h}}{\p \bar\eta^q}}, \widehat{\bar\eta^l}\right] = -\delta^l_q, \ \left[\widehat{\frac{\p\Xi_{h}}{\p z^p}}, \widehat{\eta^k}\right] = \left[\widehat{\frac{\p\Xi_{h}}{\p \bar z^q}}, \widehat{\bar\eta^l}\right] = 0. 
\]
\end{lemma}
{\it Proof.} We have
\begin{eqnarray*}
   \widehat{\frac{\p\Xi_{h}}{\p z^p}} = \frac{\p\Psi_{h}}{\p z^p} +  \frac{\p}{\p z^p} \log g
-  \left(h g^{\bar lk}\frac{\p}{\p \bar z^l}\right)\circ \frac{\p^2 \Psi_{h}}{\p z^k\p z^p}+\\
 \frac{1}{h}g_{p\bar l}h g^{\bar lk}\left(\frac{\p}{\p z^k} - \frac{\p \Psi_{h}}{\p z^k}\right) = 
\frac{\p\Psi_{h}}{\p z^p} + \frac{\p}{\p z^p}\log g  - 
h g^{\bar lk}\frac{\p^2 \Psi_{h}}{\p z^k\p z^p}\frac{\p}{\p \bar z^l} -\\
 g^{\bar lk}\frac{\p}{\p z^k}g_{pl} + \frac{\p}{\p z^p} - \frac{\p\Psi_{h}}{\p z^p}= 
 \frac{\p}{\p z^p} - h g^{\bar lk}\frac{\p^2 \Psi_{h}}{\p z^k\p z^p}\frac{\p}{\p \bar z^l}. 
\end{eqnarray*}
A similar calculation provides the formula for the operator $\widehat{\frac{\p\Xi_{h}}{\p \bar z^q}}$ 
from  Eqn.(\ref{E:hatxibar}). 
Further,
\[
   \left[\widehat{\frac{\p\Xi_{h}}{\p \eta^p}}, \widehat{\eta^k}\right] = \left[\frac{\p\Psi_{h}}{\p z^p}, - h g^{\bar lk}\frac{\p}{\p\bar z^l}\right] = g_{p\bar l}\, g^{\bar lk} = \delta^k_p.
\]
The following calculation is based on the identity
\[
       \frac{\p g^{\bar lk}}{\p z^p} = - g^{\bar lm}\frac{\p g_{m\bar n}}{\p z^p}g^{\bar nk}
\]
and Eqn. (\ref{E:kp}):
\begin{eqnarray*}
   \left[\widehat{\frac{\p\Xi_{h}}{\p z^p}}, -\frac{1}{h}\widehat{\eta^k}\right] = \left[\frac{\p}{\p z^p} - h g^{\bar lm} \frac{\p^2\Psi_{h}}{\p z^m \p z^p}\frac{\p}{\p\bar z^l}, g^{\bar nk}\frac{\p}{\p \bar z^n}
\right] = \\
\frac{\p g^{\bar nk}}{\p z^p}\frac{\p}{\p \bar z^n} -
h g^{\bar lm} \frac{\p^2\Psi_{h}}{\p z^m \p z^p}\frac{\p g^{\bar nk}}{\p\bar z^l}\frac{\p}{\p\bar z^n} + \\
h\, g^{\bar nk}\frac{\p g^{\bar lm}}{\p \bar z^n}\frac{\p^2\Psi_{h}}{\p z^m \p z^p}\frac{\p}{\p\bar z^l}+ g^{\bar nk}g^{\bar lm}\frac{\p g_{m\bar n}}{\p z^p} \frac{\p}{\p\bar z^l} = 0.
\end{eqnarray*}
The rest of the Lemma can be proved by similar calculations.

Using Lemmas \ref{L:lfrf} and  \ref{L:commrel} we will calculate the left and right multiplication operators with respect to the product $*_{h}$ for a number of symbols.

\begin{lemma}\label{L:leftast}
For a holomorphic function $a$ and an antiholomorphic function 
$b$ on $U$ the following formulas hold:
\begin{eqnarray*}
L_{a} =a,\ L_{\eta^p} = \eta^p,\ R_b = b,\ R_{\bar \eta^q} = \bar\eta^q,\\
L_{\frac{\p\Xi_{h}}{\p \eta^p}} = L_{\frac{\p\Psi_{h}}{\p z^p}}= \frac{\p}{\p \eta^p}+  \frac{\p\Psi_{h}}{\p \eta^p} = \frac{\p}{\p \eta^p}+  \frac{\p\Xi_{h}}{\p \eta^p},\\
R_{\frac{\p\Xi_{h}}{\p\bar \eta^q}} = R_{\frac{\p\Psi_{h}}{\p\bar z^q}} = \frac{\p}{\p\bar \eta^q}+  \frac{\p\Psi_{h}}{\p\bar z^q} =
\frac{\p}{\p\bar \eta^q}+  \frac{\p\Xi_{h}}{\p\bar \eta^q},\\
L_{\frac{\p\Xi_{h}}{\p z^p}} = \frac{\p}{\p z^p}+  \frac{\p\Xi_{h}}{\p z^p},\ 
R_{\frac{\p\Xi_{h}}{\p \bar z^q}} = \frac{\p}{\p\bar z^q}+  \frac{\p\Xi_{h}}{\p\bar z^q}.
\end{eqnarray*}
\end{lemma}
{\it Proof.} The four formulas in the first line are obvious. The formulas for $L_{\p\Xi_{h}/\p\eta^p}$ 
and $L_{\p\Xi_{h}/\p\bar\eta^q}$ follow from Lemma \ref{L:lfrf}. To prove the formula for $L_{\p\Xi_{h}/\p z^p}$ we will calculate first its particular case with the use of Eqns.(\ref{E:etaf}) and (\ref{E:raflbf}):
\begin{eqnarray}\label{E:partic}
   L_{\frac{\p\Xi_{h}}{\p z^p}}f = \left(\frac{\p\Psi_{h}}{\p z^p} + \frac{\p}{\p z^p}\log g + \eta^k \frac{\p^2\Psi_{h}}{\p z^k\p z^p}\right)f + \\
\frac{1}{h} g_{p\bar l} *_{h} \bar \eta^l *_{h} f = 
 \frac{\p\Xi_{h}}{\p z^p}f + \frac{1}{h} g_{p\bar l} h \, g^{\bar lk}\frac{\p f}{\p z^k} =  \frac{\p\Xi_{h}}{\p z^p}f + \frac{\p f}{\p z^p}.\nonumber
\end{eqnarray}
It follows from Eqn.(\ref{E:partic}) that
\begin{equation}\label{E:strange}
     \widehat{\frac{\p\Xi_{h}}{\p z^p}}\hat f = \widehat{\left(\frac{\p\Xi_{h}}{\p z^p}f\right)} + \widehat{\frac{\p f}{\p z^p}}. 
\end{equation}
Consider a symbol $P = u(\eta)fv(\bar \eta)$, where $u(\eta)$ and $v(\bar \eta)$ are monomials in the variables $\{\eta^k\}$ and $\{\bar \eta^l\}$, respectively. Taking into account Eqn.(\ref{E:strange}) and Lemma \ref{L:commrel}, we obtain
\begin{eqnarray*}
  \widehat{\frac{\p\Xi_{h}}{\p z^p}}\hat P = \widehat{\frac{\p\Xi_{h}}{\p z^p}}u(\hat{\eta})\hat f v(\widehat{\bar \eta}) = u(\hat{\eta})\widehat{\frac{\p\Xi_{h}}{\p z^p}}\hat f v(\widehat{\bar \eta}) = \\
u(\hat{\eta})\widehat{\left(\frac{\p\Xi_{h}}{\p z^p}f\right)}v(\widehat{\bar \eta}) +
 u(\hat{\eta})\widehat{\frac{\p f}{\p z^p}}v(\widehat{\bar \eta}) = \widehat{\left(\frac{\p\Xi_{h}}{\p z^p}P\right)} + \widehat{\frac{\p P}{\p z^p}},
\end{eqnarray*}
whence the formula for $L_{\p\Xi_{h}/\p z^p}$ follows. The formula for 
$L_{\p\Xi_{h}/\p \bar z^q}$ can be checked by similar calculations.

Our next task is to find an explicit expression for the operator $L_{\bar\eta^q}$ which is possible due to Lemmas \ref{L:lfrf} and \ref{L:leftast}.
It follows from Lemma \ref{L:leftast} that
\[
   \frac{\p \Xi_{h}}{\p z^p} = \frac{\p \Psi_{h}}{\p z^p} + \frac{\p}{\p z^p}\log g + 
\eta^k *_{h} \frac{\p^2 \Psi_{h}}{\p z^k \p z^p} + \frac{1}{h} g_{p\bar l} *_{h} \bar\eta^l.
\]
Therefore,
\begin{equation}\label{E:LLL}
     L_{\frac{\p \Xi_{h}}{\p z^p}} = L_{\left(\frac{\p \Psi_{h}}{\p z^p} + \frac{\p}{\p z^p}\log g\right)} + \eta^k L_{\frac{\p^2 \Psi_{h}}{\p z^k \p z^p}} + \frac{1}{h} L_{g_{p\bar l}}L_{\bar\eta^l}.
\end{equation}
Assuming summation over repeated indices and taking into account Eqn.(\ref{E:separ}), we may write:
\begin{equation}\label{E:delta}
    L_{g^{\bar qp}}L_{g_{p\bar l}} = L_{\left(g^{\bar qp}*_{h} g_{p\bar l}\right)} = L_{g^{\bar qp}g_{p\bar l}} = L_{\delta_l^q} =  \delta_l^q.
\end{equation}
Applying Eqn.(\ref{E:delta}) to Eqn.(\ref{E:LLL}) we obtain an expression for the operator $L_{\bar\eta^q}$:
\begin{equation}\label{E:lbaretaq}
   L_{\bar\eta^q} = h\, L_{g^{\bar qp}}\left(L_{\frac{\p \Xi_{h}}{\p z^p}} - L_{\left(\frac{\p \Psi_{h}}{\p z^p} + \frac{\p}{\p z^p}\log g\right)} - \eta^k L_{\frac{\p^2 \Psi_{h}}{\p z^k \p z^p}}\right).
\end{equation}
It can be shown that the product $*_h$ defined locally on the fibrewise polynomial functions on $TM$ is coordinate invariant and determines a global product on $\P(TM)$ even when there is no global Hermitian line bundle on $M$ with the global differential operators on it which would have these functions as their symbols. Denote by $\D(TM)\langle h\rangle$ the space of series of the form
\[
   A = \sum_{r \geq 0} h^r A_r,
\]
where $A_r \in \D(TM)$ and $\lim_{r\to\infty}A_r =0$ in the topology determined by the filtration $\{\D_r(TM)\}$. These series form an algebra. Lemma \ref{L:lfrf}, Eqns. (\ref{E:separ}) and (\ref{E:lbaretaq}) imply that for any symbol $P\in \P(TM)$ the operator $L_P$ is in the space $\D(TM)\langle h\rangle$. In particular, for $P,Q\in \P(TM)$ the product $P *_{h} Q = L_P Q$ is a polynomial in $h$, i.e., an element of $\P(TM)[h]$. Thus the space $\P(TM)[h]$ is an algebra with respect to the product $*_{h}$ and the mapping $P \mapsto L_P$ is a homomorphism from $\P(TM)[h]$ to $\D(TM)\langle h\rangle$. 

Given a local potential $\Phi_{-1}$ of the pseudo-K\"ahler form $\omega_{-1}$ on $M$, denote
\[
          \Psi = \frac{1}{\nu}\Phi_{-1}.
\]
Introduce a local potential $\Xi_{-1}$ on $TM$ as follows:
\[
  \Xi_{-1} = \Phi_{-1} + \frac{\p \Phi_{-1}}{\p z^k}\eta^k + \frac{\p\Phi_{-1}}{\p\bar z^l}\bar\eta^l
\]
and set
\begin{equation}\label{E:formpot}
   \Xi = \frac{1}{\nu}\, \Xi_{-1} + \log g = \Psi + \frac{\p\Psi}{\p z^k}\eta^k + \frac{\p\Psi}{\p\bar z^l}\bar\eta^l + \log g.
\end{equation}
A simple check shows that the form
\[
     \Omega_{-1} = -i\p_{TM}\bar\p_{TM}\Xi_{-1}
\]
is a globally defined pseudo-K\"ahler form on $TM$. The restriction of the form $\Omega_{-1}$ to the zero section $Z$ of the tangent bundle $TM$ coincides with the form $\omega_{-1}$ (under the obvious identification of the zero section $Z$ with the manifold $M$). Thus the pseudo-K\"ahler manifold $M$ is realized as a submanifold of the pseudo-K\"ahler manifold $(TM,\Omega_{-1})$ with the induced pseudo-K\"ahler structure.  It is well known that the form
\begin{equation}\label{E:omegacan}
        \omega_{\rm can} = -i\p\bar\p \log g
\end{equation} 
is globally defined on $M$. Denote by $\Omega_0$ its lift to the tangent bundle $TM$ and set
\[
   \Omega = \frac{1}{\nu}\, \Omega_{-1} + \Omega_0.
\]
There exists a deformation quantization with separation of variables on $(TM, \Omega_{-1})$ with the characterizing form $\Omega$. The corresponding star product will be denoted $*$. In the rest of the 
paper we will extend the star product $*$ on $TM$ to singular symbols supported on the zero section $Z$ which can be interpreted as a reduction to a deformation quantization with separation of variables on $M$.

Consider the `formalization' mapping $\F$ that replaces $h$ with $\nu$ (it will be used somewhat loosely). The formalizer $\F$ determines a homomorphism of the algebra $\D(TM)\langle h\rangle$ to the algebra of formal differential operators on $TM$. In particular, it maps all left and right multiplication operators with respect to the product $*_{h}$ from Lemma \ref{L:leftast} to the corresponding left and right multiplication operators with respect to the star product $*$. Thus the image of the algebra of left multiplication operators of the algebra $(\P(TM)[h],*_{h})$ with respect to $\F$ commutes with the operators $R_b,R_{\bar \eta^q}, R_{\p\Xi/\p \bar z^q},$ and $R_{\p\Xi/\p\bar\eta^q}$ and therefore belongs to the algebra of left multiplication operators with respect to the star product $*$. This implies that $\F$ determines a homomorphism of the algebra $(\P(TM)[h],*_{h})$ to the star algebra $(C^\infty(TM)[\nu^{-1},\nu]],*)$.

The pseudo-K\"ahler metric corresponding to the form $\Omega_{-1}$ has the signature $(2m,2m)$, where $m$ is the complex dimension of $M$. Since this metric is indefinite, there are no known analytic constructions of symbols on $TM$ generalizing the algebraic construction given in this section to wider classes of symbols even if $M$ is K\"ahler.

\section{Deformation quantization with separation of variables on the tangent bundle $TM$}

A number of formulas from Section \ref{S:symb} have their formal analogues (with $h$ replaced by $\nu$) which can be proved directly. We list those of them which will be used in the sequel but give no proofs.
The formal analogue of the first equation from (\ref{E:separ}) is that for $\phi,\psi \in C^\infty(M)$
\begin{equation}\label{E:recall}
     \phi * \psi = \phi\psi.
\end{equation}
The formal analogues of the other two equations from (\ref{E:separ}) follow from the definition of deformation quantization with separation of variables. Denote by $[\cdot,\cdot]_*$ the commutator with respect to the star product $*$.  The formal analogue of Eqn. (\ref{E:commut}) is as follows:
\begin{equation}\label{E:recall2}
    \left[\eta^k, \phi\right]_* = -\nu g^{\bar lk}\frac{\p\phi}{\p\bar z^l} \mbox{ and } \left[\bar\eta^l, \phi\right]_* = \nu g^{\bar lk}\frac{\p\phi}{\p z^k}. 
\end{equation}
Formulas (\ref{E:recall2}) can be written in an equivalent form:
\begin{equation}\label{E:recall3}
f * \eta^p = \eta^p f + \nu\, g^{\bar lp}\frac{\p f}{\p\bar z^l} \mbox{ and }
\bar\eta^q * f = f\bar\eta^q + \nu\, g^{\bar qk}\frac{\p f}{\p z^k}.
\end{equation}

Given a finction $f$ on $M$ (which is identified with its lift to $TM$), one can express the operators $L_f$ and $R_f$ in terms of the locally defined formal differential operators
\begin{equation}\label{E:JK}
    J = \exp\left(\nu g^{\bar lk} \frac{\p^2}{\p \eta^k \p\bar z^l}\right) \mbox{ and } K = \exp\left(\nu g^{\bar lk} \frac{\p^2}{\p z^k \p\bar \eta^l}\right)
\end{equation}
as follows,
\begin{equation}\label{E:lfjfj2}
     L_f = JfJ^{-1} \mbox{ and } R_f = KfK^{-1}.
\end{equation}

Denote by $B_*$ the formal Berezin transform of the star product $*$. Using Eqn. (\ref{E:recall}) we can prove the following lemma.
\begin{lemma}\label{L:props}
Given a local function $f = f(z,\bar z)$ and monomials $u = u(\eta), v = v(\bar\eta)$, the formulas
\[
    B_*(u f) = f * u,\  B_*(f v) = v * f, \mbox{ and } B_*(f) =f
\]
hold.
\end{lemma}
{\it Proof.}
It is sufficient to prove the lemma for a function $f$ of the form $f = ab$, where $a=a(z)$ is a holomorphic and $b=b(\bar z)$ an antiholomorphic function. It follows from Eqns. (\ref{E:ber}) and (\ref{E:recall}) that
\begin{eqnarray*}
   B_*(uf) = B_*(uab) = b * (ua) = b * (a * u) =\\
 (b * a) * u = (ab) * u = f * u.
\end{eqnarray*}
The second formula can be proved similarly. They both imply the third one. 

In the rest of this section we will calculate the canonical formal trace density of the star product $*$ on $TM$.
Given a potential (\ref{E:formpot}) of the formal form $\Omega$, we will show that the potential
\[
   \tilde\Xi = -2m \log \nu -\Psi - \frac{\p\Psi}{\p z^k}\eta^k - \frac{\p\Psi}{\p\bar z^l}\bar\eta^l + \log g
\]
satisfies Eqn. (\ref{E:normalization}) rewritten in the notations adapted to the pseudo-K\"ahler manifold $TM$ as follows:
\begin{eqnarray}\label{E:norm2}
B_*\left(\frac{\p \tilde\Xi}{\p z^p}\right)= - \frac{\p \Xi}{\p z^p},\ B_*\left(\frac{\p \tilde\Xi}{\p \eta^p}\right)=  - \frac{\p \Xi}{\p \eta^p},\ B_*\left(\frac{\p \tilde\Xi}{\p\bar z^q}\right)=  - \frac{\p \Xi}{\p\bar z^q},\nonumber\\
B_*\left(\frac{\p \tilde\Xi}{\p\bar \eta^q}\right)=  - \frac{\p \Xi}{\p\bar \eta^q}, \mbox{ and }B_*\left(\frac{d\tilde\Xi}{d\nu}\right)= - \frac{d\Xi}{d\nu}.    
\end{eqnarray}
The calculations below are based on Lemma \ref{L:props}. The following calculation proves the first formula:
\begin{eqnarray*}
B_*\left(\frac{\p \tilde\Xi}{\p z^p}\right)=
B_*\left(-\frac{\p\Psi}{\p z^p} - \eta^k \frac{\p^2 \Psi}{\p z^k \p z^p} - \frac{1}{\nu} g_{p\bar l}\bar\eta^l + \frac{\p}{\p z^k}\log g\right)=\\
-\frac{\p\Psi}{\p z^p} - \frac{\p^2 \Psi}{\p z^k \p z^p}*\eta^k  - \frac{1}{\nu}\bar \eta^l *  g_{p\bar l} + \frac{\p}{\p z^k}\log g =
-\frac{\p\Psi}{\p z^p} - \\
 \eta^k \frac{\p^2 \Psi}{\p z^k \p z^p} -
\nu g^{\bar lk}\frac{\p g_{k\bar l}}{\p z^p} - \frac{1}{\nu} g_{p\bar l}\bar\eta^l - g^{\bar lk}\frac{\p g_{p\bar l}}{\p z^k} + \frac{\p}{\p z^k}\log g = - \frac{\p \Xi}{\p z^p}.
\end{eqnarray*}
The following one proves the second formula:
\[
 B_*\left(\frac{\p \tilde\Xi}{\p \eta^p}\right)= B_*\left(-\frac{\p \Psi}{\p z^p}\right)=
-\frac{\p \Psi}{\p z^p} = - \frac{\p \Xi}{\p \eta^p}.
\]
The next two formulas can be proved similarly. The proof of the last formula in (\ref{E:norm2}) is as follows:
\begin{eqnarray*}
B_*\left(\frac{d\tilde\Xi}{d\nu}\right)=
  B_*\left(\frac{1}{\nu}\Psi + \frac{1}{\nu}\frac{\p\Psi}{\p z^k}\eta^k + \frac{1}{\nu}\frac{\p\Psi}{\p \bar z^l}\bar\eta^l - 2\frac{m}{\nu}\right)= \\
\frac{1}{\nu}\Psi +
 \frac{1}{\nu}\frac{\p\Psi}{\p z^k}*\eta^k + \frac{1}{\nu}\bar\eta^l * \frac{\p\Psi}{\p \bar z^l} - 2\frac{m}{\nu} = \frac{1}{\nu}\Psi + \\
\frac{1}{\nu}\frac{\p\Psi}{\p z^k}\eta^k + 
\frac{1}{\nu}g^{\bar lk}g_{k\bar l} + 
 \frac{1}{\nu}\frac{\p\Psi}{\p \bar z^l}\bar\eta^l + \frac{1}{\nu}g^{\bar lk}g_{k\bar l} - 2\frac{m}{\nu} = - \frac{d\Xi}{d\nu}.
\end{eqnarray*}
According to Eqn. (\ref{E:tracedensity}), the canonical trace density $\mu_*$ of the star product $*$ is given by the formula
\begin{equation}\label{E:lambda}
              \mu_* = \lambda_m e^{\left(\Xi + \tilde \Xi\right)} dz d\bar z d\eta d\bar\eta = \frac{\lambda_m}{\nu^{2m}}g^2 dz d\bar z d\eta d\bar\eta,
\end{equation}
where $\lambda_m$ is a constant and $d\eta d\bar \eta = d\eta^1\ldots d\eta^m d\bar \eta^1\ldots d\bar\eta^m$. Eventually, we obtain from Eqn. (\ref{E:mu}) that
\begin{equation}\label{E:mustar}
          \mu_* = \frac{1}{\nu^{2m}(2m)!}\Omega_{-1}^{2m},
\end{equation}
which means that the star-product $*$ is `closed' (see \cite{CFS}). The constant $\lambda_m$ is thus determined by Eqns. (\ref{E:lambda}) and (\ref{E:mustar}) and can be explicitly calculated.

\section{A fibrewise Fourier transformation}

In this section we will use the following terminology and facts. A generalized function is a functional on the smooth compactly supported densities.
If $E\to M$ is a fibre bundle then a fibrewise generalized function (i.e. a smooth family of generalized functions on the fibres of $E$) is a generalized function on the total space $E$. The action of a differential operator on the functions on $E$ extends to the space of fibrewise generalized functions on $E$ and to the space of all generalized functions on $E$. 

Denote by $\E^r$ the space of fibrewise generalized functions on the tangent bundle $TM$ supported on the zero section of $TM$ of order not greater than $r$ and set $\E = \cup_r \E^r$. Consider the subspace $\E_{\nu} \subset\E[\nu^{-1},\nu]]$ consisting of the elements such that for each $A\in \E_{\nu}$ there is an integer $s$ for which $A$ can be represented as 
\[
A = \sum_{r \geq s}\nu^r A_r 
\]
with $A_r \in \E^{r-s}$ for all $r\geq s$.

Natural formal differential operators act on the space $\E_\nu$. We will prove a more general statement.
\begin{lemma}\label{L:action}
  A formal differential operator $B$ given by the $\nu$-adically convergent series
\[
        B = \sum_{r \geq k} \nu^r B_r,
\]
where $k$ is a (possibly negative) integer and $B_r$ are natural formal differential operators, leaves invariant the space $\E_\nu$.
\end{lemma}
{\it Proof.} Given $A\in \E_\nu$, there is an integer $l$ such that $A$ can be represented as
\[
         A = \sum_{s\geq l} \nu^s A_s
\]
with $A_s \in \E^{s-l}$. Represent each natural formal differential operator $B_r$ as
\[
     B_r = \sum_{t \geq 0} \nu^t B_{r,t},
\]
where $B_{r,t}$ is a differential operator of order not greater than $t$. Now
\[
    B(A) = \sum_{r \geq k}\sum_{t \geq 0}\sum_{s \geq l} \nu^{r+ s + t} B_{r,t}(A_s).
\]
Set $n = r+s+t$. The order of the generalized function $B_{r,t}(A_s)$ is not greater than $t + s - l$. The Lemma follows from the fact that $t + s - l = n - (r + l) \leq n - (k+l)$ for any $n$.

We identify the elements of the space $C^\infty(\T,Z)$ of functions on the formal neighborhood $(\T,Z)$ of the zero section $Z$ of the cotangent bundle $\T$ with the formal series $A = A_0 + A_1 + A_2 + \ldots$, where $A_r$ is a fibrewise homogeneous polynomial of degree $r$ in the  fibre variables on $\T$ (see the Appendix for the details on formal neighborhoods). Denote by $\P^r$ the space of sums
\[
    P = \sum_{s=0}^r \frac{1}{\nu^s} P_s,
\]
where $P_s$ is a fibrewise homogeneous polynomial of degree $s$ on the cotangent bundle $\T$. The space $C^\infty(\T,Z)[\nu^{-1},\nu]]$ can be alternatively described as the set of formal series of the form
\[
       F = \sum_{r \geq s}\nu^r P_r,
\]
where $s$ is some integer and $P_r \in \P^{r-s}$ for all $r\geq s$. To see it, represent $P_r$ as
\[
          P_r = \sum_{t=0}^{r-s} \frac{1}{\nu^t} P_{r,t},
\]
where $P_{r,t}$ is a fibrewise homogeneous polynomial of degree $t$ in the fibre variables on $\T$. Then $F$ can be rewritten as an element of $C^\infty(\T,Z)[\nu^{-1},\nu]]$ as follows:
\[
       F = \sum_{r \geq s}\nu^r \left(\sum_{t=0}^{r-s} \frac{1}{\nu^t} P_{r,t}\right) = \sum_{r \geq s} \nu^r \left(\sum_{t \geq 0} P_{r+t,t}\right).
\]

The pseudo-K\"ahler metric $g_{k\bar l}$ on $M$ defines a global fibrewise density $gd\eta d\bar \eta$ on $TM$, where $g = \det(g_{k\bar l})$. Using this density and the natural pairing of $TM$ and $\T$, define a fibrewise Fourier transformation of the elements of the space $\E$ by the formula 
\begin{equation}\label{E:Four}
    \tilde A(z,\bar z, \xi,\bar \xi) = \int e^{\frac{i}{\nu}(\eta^k \xi_k + \bar \eta ^l \bar \xi_l)} A(z,\bar z, \eta,\bar \eta) g d\eta d\bar \eta.
\end{equation}
It is an isomorphism of $\E^r$ onto $\P^r$ and therefore it extends to an isomorphism of the space $\E_{\nu}$ onto $C^\infty(\T,Z)[\nu^{-1},\nu]]$. It transfers the differential operators on $TM$ to operators on  $\T$ as follows. For $f = f(z,\bar z)$ denote the pointwise multiplication operator by $f$ by the same symbol. Then
\begin{eqnarray}\label{E:etaxi}
   f \mapsto f,\ \eta^k \mapsto -i\nu \frac{\p}{\p \xi_k}, \ \bar\eta^l \mapsto -i\nu \frac{\p}{\p\bar \xi_l}, \frac{\p}{\p \eta^k} \mapsto -\frac{i}{\nu} \xi_k, \\
 \frac{\p}{\p\bar\eta^l}\mapsto -\frac{i}{\nu}\bar\xi_l, \frac{\p}{\p z^k} \mapsto \frac{\p}{\p z^k} - \frac{\p}{\p z^k}\log g, \ \frac{\p}{\p \bar z^l} \mapsto \frac{\p}{\p \bar z^l} - \frac{\p}{\p \bar z^l}\log g.\nonumber
\end{eqnarray}
Notice that the pointwise multiplication operator by a function 
$u(\eta)$ transfers to a formal differential operator
\[
     \sum_{r=0}^\infty u_r\left(-i\nu\frac{\p}{\p\xi}\right),
\]
where $\sum_{r=0}^\infty u_r$ is the Taylor series of the function $u(\eta)$ at the origin with $u_r$ a homogeneous polynomial of degree 
$r$. 

Since the star product $\ast$ on $TM$ is natural, the left and right multiplication operators $L_A$ and $R_A$  are natural for any $A \in C^\infty(TM)$. It follows from Lemma \ref{L:action} that for $A \in C^\infty(TM)[\nu^{-1},\nu]]$ the action of the operators $L_A$ and $R_A$ can be extended to the space $\E_{\nu}$. This action can be transferred to the space $C^\infty(\T,Z)[\nu^{-1},\nu]]$ via the Fourier transformation (\ref{E:Four}). Denote the corresponding operators on the space $C^\infty(\T,Z)[\nu^{-1},\nu]]$ by $\tilde L_A$ and $\tilde R_A$, respectively, so that for $B \in \E_{\nu}$
\[
   \widetilde{L_A B} = \tilde L_A \tilde B \mbox{ and } \widetilde{R_A B} = \tilde R_A \tilde B.
\]
Thus both $\E_{\nu}$  and $C^\infty(\T,Z)[\nu^{-1},\nu]]$ are bimodules over the star algebra $(C^\infty(TM)[\nu^{-1},\nu]],*)$.

It turns out quite surprisingly that the operators $\tilde L_A$ and $\tilde R_A$ are naturally expressed in terms of the so called `formal symplectic groupoid with separation of variables' over the pseudo-K\"ahler $M$ (see \cite{CMP3}). Consider the standard Poisson structure on $\T$ given on $A,B\in C^\infty(\T)$ as follows:
\[
   \{A,B\}_{\T} = \frac{\p A}{\p \xi_k}\frac{\p B}{\p z^k}-\frac{\p B}{\p \xi_k}\frac{\p A}{\p z^k} + \frac{\p A}{\p \bar\xi_l}\frac{\p B}{\p \bar z^l}-\frac{\p B}{\p \bar\xi_l}\frac{\p A}{\p \bar z^l}.
\]
There exist a Poisson and an anti-Poisson (global) morphisms $S,T: C^\infty(M) \to C^\infty(\T,Z)$, respectively, given locally by the formulas
\[
   S\phi = e^{-i\xi_k g^{\bar lk}\frac{\p}{\p \bar z^l}}\phi \mbox{ and } T\psi = e^{-i\bar\xi_l g^{\bar lk}\frac{\p}{\p z^k}}\psi.  
\]
These are the source and target mappings of the formal symplectic groupoid with separation of variables over the pseudo-K\"ahler manifold $M$. The images of the source and target mappings Poisson commute,
\[
         \{S\phi, T\psi\}_{\T} = 0.
\]
The mapping $S\otimes T: \phi\otimes\psi \mapsto (S\phi)(T\psi)$ extends to a Poisson isomorphism
\begin{equation}\label{E:isom}
    S\otimes T: C^\infty(M\times\bar M,M_{\rm diag}) \to C^\infty(\T,Z),
\end{equation}
where the factor $\bar M$ is endowed with the Poisson structure opposite to (\ref{E:PoissM}) and the product $M\times\bar M$ by the product Poisson structure (see \cite{Deq}).
For a local holomorphic function $a$ and  antiholomorphic function $b$ on $M$
\[
    Sa = a \mbox{ and } Tb = b.
\]
Given a function $f \in C^\infty(M)$ (which is identified with its lift to $\T$), its pullback via the mapping $S\otimes T$ is the function $\delta f$, where $\delta$ is the formal analytic extension mapping (see the Appendix). To check it, it is sufficient to consider a function of the form
\[
     f = \sum_i a_i b_i,
\]
where $a_i$ and $b_i$ are local holomorphic and  antiholomorphic functions, respectively.
Then
\begin{eqnarray*}
   (S\otimes T)(\delta f) = (S\otimes T)\left( \delta\sum_i a_i b_i\right) = \\
(S\otimes T)\left(\sum_i a_i \otimes b_i\right) = \sum_i S(a_i)T(b_i) = \sum_i a_i b_i = f. 
\end{eqnarray*}

Using the fact that the Fourier transforms of the operators $J$ and $K$ given by Eqn. (\ref{E:JK}) are
\begin{equation}\label{E:fourjk}
   \tilde J = \exp\left(-i \xi_k g^{\bar lk}\frac{\p}{\p\bar z^l}\right) \mbox{ and } \tilde K = \exp\left(-i \bar\xi_l g^{\bar lk}\frac{\p}{\p z^k}\right),
\end{equation}
we can prove the following lemma.
\begin{lemma}\label{L:fourst}
   Given a function $f$ on $M$, the Fourier transforms of the operators $L_f$ and $R_f$ are pointwise multiplication operators given by the formulas
\[
        \tilde L_f = S(f) \mbox{ and } \tilde R_f = T(f).
\]
\end{lemma}
{\it Proof.} Since the fibrewise Fourier transform of the pointwise multiplication operator by $f$ is also the multiplication operator by $f$, we get from Eqns. (\ref{E:lfjfj2}) and (\ref{E:fourjk}) that
\[
   \tilde L_f = \tilde J f \tilde J^{-1} = e^{-i\xi_k g^{\bar lk}\frac{\p}{\p \bar z^l}}f  = Sf.
\]
The formula for $\tilde R_f$ can be proved similarly.

Given a function $f$ on $M$, denote by $l_f$ and $r_f$ the pullbacks of the operators $\tilde L_f$ and $\tilde R_f$ via the isomorphism (\ref{E:isom}), respectively. According to Lemma \ref{L:fourst},
\begin{equation}\label{E:pullback}
       l_f = f \otimes 1 \mbox{ and } r_f = 1 \otimes f.
\end{equation}
For a Hamiltonian $A \in C^\infty(\T)$ we will denote the corresponding Hamiltonian vector field on $\T$ by $H_A$ so that for 
$B \in C^\infty(\T)$
\[
   H_A B = \{A,B\}_{\T}.
\]
Similarly, we will denote by $h_\phi$ the Hamiltonian vector field on $M$ corresponding to a Hamiltonian function $\phi$ so that 
\[
               h_\phi \psi = \{\phi,\psi\}_M
\]
for $\psi\in C^\infty(M)$. Our next task will be to calculate several left and right multiplication operators, their Fourier transforms and pullbacks via the mapping $S\otimes T$. These calculations will be used in the rest of the paper.

Since $L_{\eta^p} = \eta^p$ and $R_{\bar\eta^q} = \bar\eta^q$ we see from Eqn. (\ref{E:etaxi}) that
\begin{eqnarray}\label{E:tiltil}
  & \tilde L_{\eta^p} = -i\nu\frac{\p}{\p \xi_p} = i\nu H_{z^p} = i\nu H_{S(z^p)} = i\nu e^{-\Psi}H_{S(z^p)}e^\Psi \nonumber\\
 \mbox{ and } & \\
  & \tilde R_{\bar\eta^q} = -i\nu\frac{\p}{\p \bar\xi_q} = i\nu H_{\bar z^q} = i\nu H_{T(\bar z^q)} = i\nu e^{-\Psi}H_{T(\bar z^q)}e^\Psi.
\nonumber
\end{eqnarray}
We see from Eqn.(\ref{E:tiltil}) and the fact that $S\otimes T$ is a Poisson morphism that
\begin{eqnarray}\label{E:delta1}
    & l_{\eta^p} = i\nu \left(h_{z^p}\otimes 1\right) = i\nu e^{-\delta\Psi} \left(h_{z^p}\otimes 1\right)e^{\delta\Psi}\nonumber \\
\mbox{ and } &  \\
& r_{\bar\eta^q} = -i\nu \left(1 \otimes h_{\bar z^q}\right) = -i\nu e^{-\delta\Psi} \left(1 \otimes h_{\bar z^q}\right)e^{\delta\Psi}. \nonumber
\end{eqnarray}
Formula (\ref{E:formpot}) implies that
\begin{equation}\label{E:Lzp}
    L_{\frac{\p\Xi}{\p z^p}}=  \frac{\p\Xi}{\p z^p} + \frac{\p}{\p z^p}=     
\frac{\p\Psi}{\p z^p} + \frac{\p^2\Psi}{\p z^p \p z^k}\eta^k + \frac{1}{\nu}\, g_{p\bar l} \bar\eta^l + \frac{\p}{\p z^p}\log g + \frac{\p}{\p z^p}. 
\end{equation}
Notice that
\begin{equation}\label{E:hs}
   H_{S\left(\frac{\p\Psi}{\p z^p}\right)} = - \frac{i}{\nu}\frac{\p}{\p z^p} - \frac{\p^2\Psi}{\p z^p \p z^k} \frac{\p}{\p\xi_k} - \frac{1}{\nu}\, g_{p\bar l} \frac{\p}{\p\bar\xi_l}.
\end{equation}
It follows from Eqns.(\ref{E:etaxi}),(\ref{E:Lzp}), and (\ref{E:hs}) that
\begin{equation}\label{E:tildelfrac}
   \tilde L_{\frac{\p\Xi}{\p z^p}}= i\nu H_{S\left(\frac{\p\Psi}{\p z^p}\right)} + \frac{\p\Psi}{\p z^p} = i\nu\, e^{-\Psi} H_{S\left(\frac{\p\Psi}{\p z^p}\right)}e^\Psi. 
\end{equation}
A similar calculation shows that
\[
    \tilde R_{\frac{\p\Xi}{\p \bar z^q}}= i\nu\, e^{-\Psi} H_{T\left(\frac{\p\Psi}{\p \bar z^q}\right)}e^\Psi. 
\]
Pulling these operators back via the mapping $S\otimes T$ we obtain that
\begin{eqnarray}\label{E:delta2}
   &l_{\frac{\p\Xi}{\p z^p}}= i\nu\, e^{-\delta\Psi}\left(h_{\frac{\p\Psi}{\p z^p}} \otimes 1\right) e^{\delta\Psi} \nonumber\\
\mbox{ and }& \\
& r_{\frac{\p\Xi}{\p \bar z^q}}= -i\nu\, e^{-\delta\Psi}\left(1 \otimes h_{\frac{\p\Psi}{\p\bar z^q}}\right) e^{\delta\Psi},\nonumber
\end{eqnarray}
where $\delta$ is the formal analytic extension mapping.

Given a function $\phi$ on $M$, denote by $A(\phi)$ the following function on $TM$:
\[
               A(\phi) = \frac{\p \phi}{\p z^k}\eta^k + \frac{\p \phi}{\p \bar z^l}\bar\eta^l.
\]
Notice that $A(z^k) = \eta^k$ and $A(\bar z^l) = \bar\eta^l$. There is a formal analogue (with $h$ replaced by $\nu$) of long and indirect formula (\ref{E:lbaretaq}) for the operator $L_{\bar\eta^q}$. It turns out that the Fourier transform of this operator has a nice expression. We will find a general formula for the operator $\tilde L_{A(\phi)}$ by working first with its pullback $l_{A(\phi)}$ via the mapping $S\otimes T$.

Using Eqns. (\ref{E:PoissM}), (\ref{E:recall}), and (\ref{E:recall2}) we get that
\begin{eqnarray}\label{E:aphipsi}
   [A(\phi),\psi]_* = \left[\eta^k * \frac{\p \phi}{\p z^k} + \frac{\p \phi}{\p \bar z^l} * \bar\eta^l, \psi\right]_* = \\
\left[\eta^k, \psi\right]_* * \frac{\p \phi}{\p z^k} + \frac{\p \phi}{\p \bar z^l} *\left[\bar\eta^l,\psi\right]_* = i\nu\{\phi,\psi\}_M = i\nu h_\phi\psi.  \nonumber
\end{eqnarray}
Using Eqn. (\ref{E:pullback}), we obtain from Eqn. (\ref{E:aphipsi}) that
\begin{equation}\label{E:laotimes}
    \left[l_{A(\phi)}, \psi\otimes 1\right] = i\nu h_\phi\psi \otimes 1.
\end{equation}
Since the left multiplication operators commute with the right multiplication operators, we see from Eqns. (\ref{E:pullback}) and (\ref{E:laotimes}) that the operator
\[
    B(\phi) = l_{A(\phi)} - i\nu e^{-\delta\Psi}\left(h_\phi \otimes 1\right)e^{\delta\Psi}
\]
commutes with the functions $\psi\otimes 1$ and $1\otimes\psi$ for any $\psi\in C^\infty(M)$ and thus with the pointwise multiplication operators by the elements of the space $C^\infty(M\times \bar M,M_{\rm diag})$. Thus $B(\phi)$ is a pointwise multiplication operator itself.
It follows from Eqns. (\ref{E:delta1}) and (\ref{E:delta2}) that the multiplication operator $B(\phi)$ commutes with the operators
\[
  r_{\bar\eta^q} = -i\nu e^{-\delta\Psi} \left(1 \otimes h_{\bar z^q}\right)e^{\delta\Psi} \mbox{ and }
r_{\frac{\p\Xi}{\p \bar z^q}}= -i\nu e^{-\delta\Psi}\left(1 \otimes h_{\frac{\p\Psi}{\p\bar z^q}}\right) e^{\delta\Psi},
\]
and therefore with the operators
\[
  1 \otimes h_{\bar z^q} \mbox{ and } 1 \otimes h_{\frac{\p\Psi}{\p\bar z^q}}.
\]
This implies that the function $B(\phi)$ is of the form $C(\phi) \otimes 1$. To identify the mapping $\phi \mapsto C(\phi)$ we will push forward the formula
\[
          l_{A(\phi)} = i\nu e^{-\delta\Psi}\left(h_\phi \otimes 1\right)e^{\delta\Psi} + C(\phi) \otimes 1
\]
via the mapping $S\otimes T$ obtaining that
\begin{equation}\label{E:tildelaphi}
       \tilde L_{A(\phi)} = i\nu e^{-\Psi}\left(H_{S(\phi)}\right)e^{\Psi} + S(C(\phi)).
\end{equation}
The following formula can be obtained by straightforward calculations with the use of Eqns. (\ref{E:recall}) and (\ref{E:recall2}). For $\phi,\psi \in C^\infty(M)$
\begin{equation}\label{E:phistara}
   A(\phi\psi) = \phi * A(\psi) + \psi * A(\phi) - \nu g^{\bar lk}\left(\frac{\p\phi}{\p z^k}\frac{\p\psi}{\p\bar z^l} + \frac{\p\psi}{\p z^k}\frac{\p\phi}{\p \bar z^l}\right).
\end{equation}
Calculating $\tilde L_{A(\phi\psi)}$ in two different ways using Lemma \ref{L:fourst}, Eqns. (\ref{E:tildelaphi}) and (\ref{E:phistara}) we get that
\begin{eqnarray}\label{E:inuepsi}
i\nu e^{-\Psi}\left(H_{S(\phi\psi)}\right)e^{\Psi} + S(C(\phi\psi))
    = i\nu S(\phi)e^{-\Psi}\left(H_{S(\psi)}\right)e^{\Psi} + \nonumber\\
i\nu S(\psi)e^{-\Psi}\left(H_{S(\phi)}\right)e^{\Psi} - S\left(\nu g^{\bar lk}\left(\frac{\p\phi}{\p z^k}\frac{\p\psi}{\p\bar z^l} + 
\frac{\p\psi}{\p z^k}\frac{\p\phi}{\p \bar z^l}\right) \right). 
\end{eqnarray}
Simplifying Eqn. (\ref{E:inuepsi}) we arrive at the statement that the Hochschild differential of the operator $C$ is
\[
    d_{\rm Hoch}C(\phi,\psi) = \phi C(\psi) - C(\phi\psi) + C(\phi)\psi = \nu g^{\bar lk}\left(\frac{\p\phi}{\p z^k}\frac{\p\psi}{\p\bar z^l} + 
\frac{\p\psi}{\p z^k}\frac{\p\phi}{\p \bar z^l}\right)
\]
and thus coincides with the Hochschild differential of the operator $-\nu\Delta$, where $\Delta$ is the Laplace-Beltrami operator (\ref{E:Delta}) on $M$. Therefore
\[
    C = -\nu\Delta + D,
\] 
where $D$ is a derivation. To determine $D$, we use Eqns. (\ref{E:tiltil}) and (\ref{E:tildelfrac}). We see from Eqns. (\ref{E:tiltil}) and (\ref{E:tildelaphi}) that $Dz^k=0$, i.e. that $D$ differentiates in antiholomorphic directions only.
It follows from Eqn. (\ref{E:recall2}) that
\[
   \left[\eta^k, \frac{\p^2 \Psi}{\p z^k\p z^p} \right]_* = - \frac{\p}{\p z^p}\log g, 
\]
whence 
\begin{equation}\label{E:onehand}
   \frac{\p\Xi}{\p z^p} = \frac{\p\Psi}{\p z^p} + \frac{\p^2 \Psi}{\p z^k\p z^p} * \eta^k + \frac{1}{\nu}g_{p\bar l}*\bar\eta^l.
\end{equation}
We see from  Eqns. (\ref{E:tiltil}), (\ref{E:tildelfrac}), and (\ref{E:onehand}) that
\begin{eqnarray}\label{E:longeqlty}
   i\nu\, e^{-\Psi} H_{S\left(\frac{\p\Psi}{\p z^p}\right)}e^\Psi = S\left(\frac{\p\Psi}{\p z^p}\right) + i\nu S\left(\frac{\p^2 \Psi}{\p z^k\p z^p}\right)e^{-\Psi}H_{S\left(z^p\right)}e^\Psi +\nonumber\\
 \frac{1}{\nu}
S\left(g_{p\bar l}\right)\left(i\nu e^{-\Psi}H_{S\left(\bar z^l\right)}e^\Psi + S\left(D(\bar z^l)\right)\right). 
\end{eqnarray}
Simplifying Eqn. (\ref{E:longeqlty}) we get that
\[
          D(\bar z^l) = -\nu g^{\bar lk} \frac{\p\Psi}{\p z^k},
\]
which means that
\[
        D = -\nu g^{\bar lk} \frac{\p\Psi}{\p z^k}\frac{\p}{\p \bar z^l}.
\]
Eventually we obtain a formula for the operator $\tilde L_{A(\phi)}$:
\begin{equation}\label{E:lfinal}
    \tilde L_{A(\phi)} = i\nu e^{-\Psi}\left(H_{S(\phi)}\right)e^{\Psi} - \nu S\left(\Delta\phi + g^{\bar lk} \frac{\p\Psi}{\p z^k}\frac{\p\phi}{\p \bar z^l}\right).
\end{equation}
Similarly, one can obtain a formula for the operator $\tilde R_{A(\phi)}$:
\begin{equation}\label{E:rfinal}
    \tilde R_{A(\phi)} = i\nu e^{-\Psi}\left(H_{T(\phi)}\right)e^{\Psi} - \nu T\left(\Delta\phi + g^{\bar lk} \frac{\p\Psi}{\p \bar z^l}\frac{\p\phi}{\p z^k}\right).
\end{equation}

\section{A product on the singular symbols}

We want to show that there exists a natural construction of an associative product $\bullet$ on the space $C^\infty(\T,Z)[\nu^{-1},\nu]]$ compatible with the bimodule structure over the star algebra $(C^\infty(TM)[\nu^{-1},\nu]],*)$ so that one can define an associative product on the direct sum 
\[
  C^\infty(TM)[\nu^{-1},\nu]]\oplus C^\infty(\T,Z)[\nu^{-1},\nu]]. 
\]
The compatibility conditions are as follows. For $F \in C^\infty(TM)[\nu^{-1},\nu]]$ and $A,B \in C^\infty(\T,Z)[\nu^{-1},\nu]]$,     
\begin{eqnarray}\label{E:comp}
   (\tilde L_F A) \bullet B = \tilde L_F(A\bullet B), (\tilde R_F A) \bullet B = A \bullet (\tilde L_F B), \mbox{ and }\\
 A \bullet (\tilde R_F B) = \tilde R_F(A\bullet B).\nonumber
\end{eqnarray}
Assume now that there exists such an associative product $\bullet$ on the space $C^\infty(\T,Z)[\nu^{-1},\nu]]$ satisfying the compatibility conditions (\ref{E:comp}). 
We want to study the properties that this product must have.
We identify the functions on $\T$ which do not depend on the fibre variables $\xi_k, \bar\xi_l$ with the functions on $M$. Thus the space $C^\infty(M)[\nu^{-1},\nu]]$ can be treated as a subspace of $ C^\infty(\T,Z)[\nu^{-1},\nu]]$.

Since $\tilde L_{\eta^k} = - i\nu\frac{\p}{\p \xi_k}$ and $\tilde R_{\bar\eta^l} = - i\nu\frac{\p}{\p \bar\xi_l}$, then it follows from Eqn. (\ref{E:comp}) that
for $\phi \in C^\infty(M)[\nu^{-1},\nu]]$  and $A\in C^\infty(\T,Z)[\nu^{-1},\nu]]$  the product $\phi \bullet A$ does not depend on the holomorphic fibre variables $\xi_k$ and $A \bullet \phi$ does not depend on the antiholomorphic fibre variables $\bar\xi_l$. Therefore $C^\infty(M)[\nu^{-1},\nu]]$  is closed with respect to the product $\bullet$. We will make an assumption that the algebra $(C^\infty(M)[\nu^{-1},\nu]], \bullet)$ has a unity $\e$ which is an invertible formal function,
\[
   \e = \sum_{r \geq n}\nu^r \e_r, 
\] 
where $n$ is an integer and $\e_n$ is nonvanishing. In the whole algebra $(C^\infty(\T,Z)[\nu^{-1},\nu]],\bullet)$ the element $\e$ will be an idempotent.
Taking into account Eqn. (\ref{E:comp}) we see that for $f \in C^\infty(M)[\nu^{-1},\nu]]$ 
the element $\e \bullet (S(f)\e) = (T(f)\e) \bullet \e$ also belongs to the space $C^\infty(M)[\nu^{-1},\nu]]$. Define a linear operator $Q$ on the space $C^\infty(M)[\nu^{-1},\nu]]$ by the formula
\begin{equation}\label{E:defQ}
    Q(f)\e = \e \bullet (S(f)\e) = (T(f)\e) \bullet \e.
\end{equation}
Notice that $Q(1) = 1$. We will assume that $Q$ is a formal differential operator. Define an associative product $\circ$ on $C^\infty(M)[\nu^{-1},\nu]]$ by the following formula. For $\phi,\psi \in C^\infty(M)[\nu^{-1},\nu]]$ set
\[
     (\phi\circ\psi)\e = (\phi\e) \bullet  (\psi\e).
\] 
The unit constant $1$ is the unity in the algebra $(C^\infty(M)[\nu^{-1},\nu]],\circ)$. We will make a further assumption that the operation $\circ$ is a star product on $M$ with respect to some Poisson structure on $M$. The assumptions we have made allow us to identify this star product. Since for a local holomorphic function $a$ and an antiholomorphic function $b$ 
\[
    \tilde L_a = S(a) = a \mbox{ and } \tilde R_b = T(b) = b,
\]
we see from Eqn. (\ref{E:comp}) that
\[
   (a\circ \psi)\e = (a\e) \bullet (\psi\e) = a(\e\bullet(\psi\e)) = a\psi\e,
\]
whence $a \circ \psi = a\psi$ and, similarly, $\phi \circ b = \phi b$. This means that the star product $\circ$ defines a deformation quantization with separation of variables with respect to some K\"ahler-Poisson structure on $M$ given by a bivector field $\eta$ of type $(1,1)$ with respect to the complex structure on $M$. 

Denote by $B_\circ$ the formal Berezin transform corresponding to the star product $\circ$ so that $B_\circ(ab) = b \circ a$, where $a,b$ are as above. Introduce an equivalent star product $\circ'$ on $M$ by the formula
\begin{equation}\label{E:circ1}
     \phi \circ' \psi = B_\circ^{-1} (B_\circ(\phi)\circ B_\circ(\psi)).
\end{equation}
This is the star product of a deformation quantization with separation of variables on the K\"ahler-Poisson manifold $(\bar M,\eta)$ with the opposite complex structure, so that locally $b \circ' f = bf$ and $f \circ' a = af$. Denote by $\tilde\circ$ the opposite star product, $\tilde\circ = (\circ')^{\rm opp}$, so that
\begin{equation}\label{E:circ2}
     \phi\, \tilde\circ\, \psi = B_\circ^{-1} (B_\circ(\psi)\circ B_\circ(\phi)).
\end{equation} 
This is the star product of a deformation quantization with separation of variables on $(M,-\eta)$. The star products $\circ$ and $\tilde\circ$ are dual.

It follows from the definition (\ref{E:defQ}) of the operator $Q$ and compatibility conditions (\ref{E:comp}) that
\begin{eqnarray*}
   Q(ab)\e = \e\bullet (S(ab)\e) = (T(b)\e) \bullet (S(a)\e) = \\
(b\e)\bullet(a\e) = (b \circ a)\e = B_\circ(ab)\e.
\end{eqnarray*}
Therefore $Q=B_\circ$ and we obtain the following formula
\begin{equation}\label{E:Berez}
     \e \bullet (S(f)\e) = (T(f)\e) \bullet \e = B_\circ(f)\e,
\end{equation}
where $f\in C^\infty(M)[\nu^{-1},\nu]]$. Introduce a formal form
\begin{equation}\label{E:dual}
    \tilde\omega = -\frac{1}{\nu}\, \omega_{-1} + \omega_{\rm can},
\end{equation}
where $\omega_{\rm can}$ is given by Eqn. (\ref{E:omegacan}).
Denote by $\tilde\star$ the star product of the deformation quantization with separation of variables on the pseudo-K\"ahler manifold $(M, -\omega_{-1})$ whose characterizing form is $\tilde\omega$ and by $\star'$ the opposite star product, $\star' = (\tilde \star)^{\rm opp}$. Let $\star$ denote the star product dual to $\tilde\star$. It is a deformation quantization with separation of variables on the 
pseudo-K\"ahler manifold $(M, \omega_{-1})$ whose characterizing form will be denoted $\omega$. Our goal is to show that the star products $\circ$ and $\star$ must coincide.

Fix a contractible coordinate chart $U \subset M$ with holomorphic coordinates $\{z^k,\bar z^l\}$. Denote by $\Phi_{-1}$ a potential of the form $\omega_{-1}$ on $U$ so that $\omega_{-1} = -i\p\bar\p \Phi_{-1}$ and set
\[
   \Psi = (1/\nu)\Phi_{-1}.
\]
The formal invertible function $\e$ can be represented on $U$ as $\e = e^\theta$ for some formal function $\theta = n \log \nu + \theta_0 + \nu \theta_1 +\ldots$. Set 
\begin{equation}\label{E:Phi}
\Phi = \Psi + \theta = (1/\nu)\Phi_{-1} + n\log \nu + \theta_0 + \nu\theta_1 + \ldots
\end{equation}
 so that $\e = e^{\Phi-\Psi}$. 
We see from Eqn. (\ref{E:lfinal}) that 
\begin{equation}\label{E:leftoper}
\tilde L_{\bar\eta^l} = i\nu e^{-\Psi} H_{S(\bar z^l)}e^\Psi - S\left(\nu g^{\bar lk}\frac{\p \Psi}{\p z^k}\right).
\end{equation}

We will need the following two technical lemmas.
\begin{lemma}\label{L:techn}
Given a function $f=f(z,\bar z)$ on the chart $U$, the following formulas hold:
\begin{eqnarray*}
H_{S(\bar z^l)}S(f) = \{S(\bar z^l),S(f)\}_{\T} = S\left(-i g^{\bar l k}\frac {\p f}{\p z^k}\right);\\
   H_{S(\bar z^l)}f=\{S(\bar z^l),f\}_{\T} =  S(-i g^{\bar l k})\frac {\p f}{\p z^k}.
\end{eqnarray*}
\end{lemma}
{\it Proof.}
\[
   \{S(\bar z^l),S(f)\}_{\T} = S(\{\bar z^l,f\}_M) = S\left(-i g^{\bar l k}\frac {\p f}{\p z^k}\right).
\]
Using the first formula, we obtain
\begin{eqnarray*}
\{S(\bar z^l), f\}_{\T} = \frac{\p S(\bar z^l)}{\p \xi_k} 
\frac {\p f}{\p z^k} = \{S(\bar z^l), z^k\}_{\T}\frac {\p f}{\p z^k} = \\
 \{S(\bar z^l), S(z^k)\}_{\T}\frac {\p f}{\p z^k} = 
S(-i g^{\bar l k})\frac {\p f}{\p z^k},
\end{eqnarray*}
which concludes the proof.
\begin{lemma}\label{L:techn2}
 Given local functions $\phi,\psi$ on $M$, the following formula holds:
\[
    \e \bullet \left(S(\phi)\psi \e\right) = \left(B_\circ (\phi) \circ \psi\right)\e = B_\circ\left( B_\circ^{-1} (\phi) \tilde\circ \psi\right)\e.
\]
\end{lemma}
{\it Proof.} It is sufficient to consider the function $\psi$ of the form
\[
     \psi = \sum_i a_i b_i,
\]
where $a_i$ and $b_i$ are local holomorphic and  antiholomorphic functions, respectively. Then, using Eqns. (\ref{E:circ1}) and (\ref{E:circ2}), we get that
\begin{eqnarray*}
\e\bullet \left(S(\phi)\psi\e\right)= \e\bullet \sum_i(S(\phi a_i)T(b_i)\e)= \\
\sum_i \left(\e\bullet S(\phi a_i)\e\right)T(b_i)= \sum_i B_\circ (\phi a_i)b_i\e = \\
\sum_i \left(B_\circ (\phi \circ' a_i)\circ b_i\right)\e = 
\left(B_\circ(\phi)\circ \sum_i(a_i\circ b_i)\right)\e = \\
\left(B_\circ(\phi)\circ \psi\right)\e = 
B_\circ \left(\phi \circ' B_\circ^{-1}(\psi)\right)\e =
B_\circ\left(B_\circ^{-1}(\psi)\tilde\circ \phi\right)\e,
\end{eqnarray*}
which proves the Lemma.

Since $\tilde R_{\bar\eta^l}= -i\nu\p/\p\bar\xi_l$ and the function $\e$ does not depend on the fibre variables $\xi,\bar\xi$, we obtain from formulas (\ref{E:Berez}), (\ref{E:leftoper}), and Lemma \ref{L:techn} that for an arbitrary formal function $f \in C^\infty(U)[\nu^{-1},\nu]]$
\begin{eqnarray}\label{E:long}
0 = (\tilde R_{\bar\eta^l}\e) \bullet S\left(\frac{1}{\nu}\, g_{p\bar l}f\right)\e = \e\bullet \tilde L_{\bar\eta^l}\left(S\left(\frac{1}{\nu}\, g_{p\bar l}f\right)\e\right) = \nonumber\\
\e\bullet i\nu e^{-\Psi} H_{S(\bar z^l)}
\left(S\left(\frac{1}{\nu}\, g_{p\bar l}f\right)e^\Phi\right) - \e\bullet S\left(g^{\bar lk}\frac{\p \Psi}{\p z^k}g_{p\bar l}f\right)=\nonumber\\
\e\bullet i \left\{S(\bar z^l),S\left(g_{p\bar l}f\right)\right\}_{\T}\e +
\e\bullet\left(iS\left(g_{p\bar l}f\right)\{S(\bar z^l),\Phi\}_{\T}\e\right) -\nonumber\\
 \e\bullet S\left(\frac{\p \Psi}{\p z^p}f\right)\e = 
\e\bullet S\left(g^{\bar lk}\frac{\p}{\p z^k}(g_{p\bar l}f)\right)\e + \\
\e\bullet \left(S\left(g_{p\bar l}f\right)S\left(g^{\bar lk}\right)\frac{\p\Phi}{\p z^k}\e\right) - 
B_\circ\left(\frac{\p\Psi}{\p z^p}f\right)\e = \nonumber\\
B_\circ\left(\frac{\p f}{\p z^p} + f\frac{\p}{\p z^p}\log g\right)\e + \e\bullet \left(S(f)\frac{\p\Phi}{\p z^p}\e\right) - B_\circ\left(\frac{\p\Psi}{\p z^p}f\right)\e.\nonumber
\end{eqnarray}
We conclude from Eqn (\ref{E:long}) that
\begin{equation}\label{E:interm}
   B_\circ\left(\frac{\p f}{\p z^p} + f\frac{\p}{\p z^p}(-\Psi + \log g)\right)\e = -\e\bullet \left(S(f)\frac{\p\Phi}{\p z^p}\e\right).
\end{equation}
Using Lemma \ref{L:techn2} we obtain that
\begin{equation}\label{E:rhs}
   \e\bullet \left(S(f)\frac{\p\Phi}{\p z^p}\e\right)=
B_\circ\left(B_\circ^{-1}\left(\frac {\p \Phi}{\p z^p}\right)\tilde\circ f\right)\e.
\end{equation}
Formulas (\ref{E:interm}) and (\ref{E:rhs}) imply that for any formal function $f$ on $U$
\begin{equation}\label{E:ltilde}
   \left(-B_\circ^{-1}\left(\frac {\p \Phi}{\p z^p}\right)\right)\tilde\circ f = \frac{\p f}{\p z^p} + f\frac{\p}{\p z^p}(-\Psi + \log g).
\end{equation}
A calculation similar to (\ref{E:long}) that starts with the observation that
\[
    0 = T\left(\frac{1}{\nu}\, g_{k\bar q}f\right)\e\bullet (\tilde L_{\eta^k}\e) = \tilde R_{\eta^k}\left(T\left(\frac{1}{\nu}\, g_{k\bar q}f\right)\e\right)\bullet \e 
\]
shows that
\begin{equation}\label{E:rtilde}
   f \tilde\circ \left(-B_\circ^{-1}\left(\frac {\p \Phi}{\p \bar z^q}\right)\right) = \frac{\p f}{\p \bar z^q} + f\frac{\p}{\p \bar z^q}(-\Psi + \log g).
\end{equation}
Since $-\Psi + \log g$ is a potential of the form $\tilde\omega$, 
it immediately follows from Eqn. (\ref{E:ltilde}) or (\ref{E:rtilde}) and the description of the star products with separation of variables on a pseudo-K\"ahler manifold that the star product $\tilde\circ$ must coincide with $\tilde\star$ and thus the star products $\star$ and $\circ$ must coincide as well. Denote by $B_{\tilde\circ},B_\star$, and $B_{\tilde\star}$ the formal Berezin transforms of the star products $\tilde\circ,\star$, and $\tilde\star$, respectively. Thus we must have that 
$B_\star = B_\circ = B^{-1}_{\tilde\circ}= B^{-1}_{\tilde\star}$. Setting $f=1$ in Eqns. (\ref{E:ltilde}) and (\ref{E:rtilde}) and replacing $B_\circ^{-1}$ with $B_{\tilde\star}$ we get that
\begin{eqnarray}\label{E:density}
&  B_{\tilde\star}\left(\frac {\p \Phi}{\p z^p}\right) = - \frac{\p}{\p z^p}(-\Psi + \log g)
\nonumber\\
\mbox{ and }&\\ 
& B_{\tilde\star}\left(\frac {\p \Phi}{\p \bar z^q}\right) = - \frac{\p}{\p \bar z^q}(-\Psi + \log g).
\nonumber
\end{eqnarray}
Formulas (\ref{E:density}) imply that $\Phi$ must be a potential of the form $\omega$ and the density $e^{\Phi-\Psi}gdzd\bar z$ must be a local trace density of the star product $\star$. 
There exists a constant $\kappa_m$ such that
\begin{equation}\label{E:kappa}
      \frac{1}{m!}\omega_{-1}^m = \kappa_m g dzd\bar z.    
\end{equation}
Since $\e = e^{\Phi-\Psi}$, we conclude from (\ref{E:kappa}) that $\e\omega_{-1}^m$ must be a global trace density of the star product $\star$, which determines $\e$ up to a formal constant factor. We see that the assumptions made in this section determine what the product $\circ$ and the formal function $\e$ might be. 

Now we will give an explicit definition of the product $\bullet$. Denote by $\mu_\star$ the canonical formal trace density of the star product $\star$ and fix an arbitrary nonzero formal constant $C(\nu)$. There exists a unique invertible formal function
\begin{equation}\label{E:extend}
   \e = \sum_{r \geq n}\nu^r \e_r 
\end{equation}
on $M$ for some integer $n$ and with $\e_n$ nonvanishing such that
\begin{equation}\label{E:defeps}
    \mu_\star = C(\nu)\e\omega_{-1}^m.
\end{equation}
We want to define an operation $\bullet$ on $C^\infty(\T,Z)[\nu^{-1},\nu]]$ such that
\begin{equation}\label{E:Berez2}
     \e \bullet (S(f)\e) = (T(f)\e) \bullet \e = B_\star(f)\e,
\end{equation}
as suggested by Eqn. (\ref{E:Berez}).  Formulas (\ref{E:comp}) and (\ref{E:Berez2}) allow to define the product $\bullet$ on the elements of $ C^\infty(\T,Z)[\nu^{-1},\nu]]$ of the form
$S(\phi)T(\psi)\e$ with $\phi,\psi \in C^\infty(M)[\nu^{-1},\nu]]$ as follows:
\begin{equation}\label{E:bullet}
   S(\phi_1)T(\psi_1)\e \bullet S(\phi_2)T(\psi_2)\e = S(\phi_1)B_\star(\psi_1\phi_2)T(\psi_2)\e.    
\end{equation}
Using the method explained in the Appendix one can extend the product $\bullet$ to the whole space $C^\infty(\T,Z)[\nu^{-1},\nu]]$. One can also show applying the technique used in this paper that the product $\bullet$ satisfies compatibility conditions (\ref{E:comp}). Now we will prove the associativity of this product.
\begin{lemma}
  The product $\bullet$ is associative.
\end{lemma}
{\it Proof.} Because of the compatibility conditions (\ref{E:comp})  it is sufficient to prove that for any functions $\psi_1, \phi_2,\psi_2,\phi_3 \in C^\infty(M)$
\[
      \left(T(\psi_1)\e \bullet S(\phi_2)T(\psi_2)\e\right) \bullet S(\phi_3)\e = T(\psi_1) \bullet 
\left(S(\phi_2)T(\psi_2)\e \bullet S(\phi_3)\e\right)
\]
or, equivalently, that
\begin{equation}\label{E:assoc}
    B_\star(\psi_1\phi_2)T(\psi_2)\e \bullet S(\phi_3)\e = T(\psi_1)\e \bullet S(\phi_2) B_\star(\psi_2\phi_3)\e. 
\end{equation}
Using Eqn. (\ref{E:comp}), simplify Eqn. (\ref{E:assoc}) as follows:
\begin{equation}\label{E:assoc2}
    B_\star(\psi_1\phi_2)T(\psi_2\phi_3)\e \bullet \e = \e \bullet S(\psi_1\phi_2) B_\star(\psi_2\phi_3)\e.
\end{equation}
Setting $\phi = \psi_1\phi_2$ and $\psi = \psi_2\phi_3$ rewrite  Eqn. (\ref{E:assoc2}) as
\begin{equation}\label{E:assoc3}
    B_\star(\phi)T(\psi)\e \bullet \e = \e \bullet S(\phi) B_\star(\psi)\e.   
\end{equation}
It is sufficient to prove Eqn. (\ref{E:assoc3}) for $\phi = a_1\star' b_1$ and $\psi = a_2\star' b_2$, where $a_1,a_2$ are local holomorphic and $b_1,b_2$ local antiholomorphic functions on $M$. The left-hand side of Eqn. (\ref{E:assoc3}) takes the form
\begin{eqnarray*}
    a_1 b_1 T(\psi)\e \bullet \e = S(a_1)T(b_1 \psi)\e \bullet \e = S(a_1)B_\star(b_1\psi)\e =\\ 
    B_\star(a_1) \star B_\star(b_1 \star' \psi)\e = B_\star(a_1 \star' b_1 \star' \psi) \e = B_\star(\phi\star'\psi)\e. 
\end{eqnarray*}
A similar calculation shows that the right-hand side of Eqn. (\ref{E:assoc3}) also equals $B_\star(\phi\star'\psi)\e$, which proves the Lemma.

For any formal function $f \in C^\infty(M)[\nu^{-1},\nu]]$ define an element ${\bf Q}_f\in C^\infty(M)[\nu^{-1},\nu]]$ by the formula
\[
        {\bf Q}_f = f\e.
\]
\begin{theorem}\label{T:covar}
  The mapping $f \mapsto {\bf Q}_f$ is an isomorphism of the algebra $(C^\infty(M)[\nu^{-1},\nu]],\star)$ onto the algebra $(C^\infty(M)[\nu^{-1},\nu]],\bullet)$.
\end{theorem}
{\it Proof.} For functions $\phi,\psi \in C^\infty(M)[\nu^{-1},\nu]]$ we have to prove that
\begin{equation}\label{E:homom}
    {\bf Q}_\phi \bullet {\bf Q}_\psi = {\bf Q}_{\phi\star\psi}.
\end{equation}
It is sufficient to prove Eqn. (\ref{E:homom}) locally for functions $\phi,\psi$ of the form
$\phi = a_1b_1$ and $\psi = a_2b_2$, where $a_1,a_2$ are local holomorphic and $b_1,b_2$ local antiholomorphic functions. It follows from Eqns. (\ref{E:comp}) and (\ref{E:Berez2}) that
\begin{eqnarray*}
{\bf Q}_\phi \bullet {\bf Q}_\psi =  (\phi\e)\bullet(\psi\e) = (a_1b_1\e)\bullet(a_2b_2\e) =
 \left(S(a_1)T(b_1)\e\right)\bullet\\
 \left(S(a_2)T(b_2)\e\right)= S(a_1)\left(\e\bullet S(b_1a_2)\e\right)T(b_2) = 
a_1\left(B_\star(b_1a_2)\e\right)b_2 =\\
 \left(a_1 \star b_1 \star a_2 \star b_2\right)\e = 
\left((a_1b_1) \star (a_2b_2)\right)\e = 
(\phi\star\psi)\e = {\bf Q}_{\phi\star\psi},
\end{eqnarray*}
which concludes the proof.

Setting $f=1$ in Eqn.(\ref{E:Berez2}) we see that the element $\e$ is an idempotent in the algebra $(C^\infty(\T,Z)[\nu^{-1},\nu]],\bullet)$,
\[
            \e \bullet \e = \e.
\]
For any function $f \in C^\infty(M)[\nu^{-1},\nu]]$ we call the element of the space $C^\infty(\T,Z)[\nu^{-1},\nu]]$ given by Eqn. (\ref{E:Berez2}) the Toeplitz element corresponding to the function $f$ and denote it ${\bf T}_f$. Thus
\[
             {\bf T}_f = \e \bullet \left(\tilde L_f \e\right) = \left(\tilde R_f \e\right)\bullet \e,
\]
which is analogous to the definition of a Toeplitz operator. The Toeplitz elements in $C^\infty(\T,Z)[\nu^{-1},\nu]]$ are exactly the elements which do not depend on the fibre variables $\xi,\bar\xi$ and thus can be identified with the elements of $C^\infty(M)[\nu^{-1},\nu]]$. 

{\it Remark.} Here we would like to give more heuristic arguments to corroborate the analogy between the Toeplitz operators on the sections of a quantum line bundle over $M$ and the Toeplitz elements in the algebra $(C^\infty(\T,Z)[\nu^{-1},\nu]],\bullet)$. Assume that $(M,\omega_{-1})$ is a compact K\"ahler manifold and $L$ is a global quantum line bundle. The Hilbert structure on the sections of the $N$-th tensor power $L^{\otimes N}$ of $L$ is given by the norm $||\cdot ||_h$ such that
\[
      ||s||_h^2 = \int |s|_h^2\, \omega_{-1}^m,
\]
where $s$ is a section and $|\cdot|_h$ is the Hermitian fibre metric (it is implied that $h = 1/N$). The symbol mapping $P \mapsto \hat P$ constructed in Section \ref{S:symb} is involutive. Namely, the complex conjugate symbol $\bar P$ corresponds to the Hermitian conjugate operator $\hat P^*$. 
Let $\alpha = f_k(z,\bar z) dz^k$ be a global differential form of type $(1,0)$ on $M$. Then $f_k\eta^k$ is a global function on $TM$. Consider the global symbol
\[
    P = J_h (f_k\eta^k) = f_k\eta^k + h g^{\bar lk}\frac{\p f_k}{\p\bar z^l} = f_k *_h \eta^k
\]
on $TM$. The corresponding global differential operator 
\[
    \hat P = f_k\nabla^k = \nabla^\alpha
\]
annihilates the holomorphic sections of $L^{\otimes N}$. The range of the conjugate operator $\hat P^*$ with the symbol $\bar P = \bar\eta^l *_h \bar f_l$, given by the formula
\[
           \hat P^* = \nabla^{\bar l} \circ \bar f_l,
\]
is orthogonal to the space of holomorphic sections of $L^{\otimes N}$ (here $\circ$ denotes composition of operators). Thus, for any Toeplitz operator $T_\phi^{(N)}$,
\[
     \hat PT_\phi^{(N)} = T_\phi^{(N)}\hat P^* =0.
\]
This statement has an obvious analogue for the Toeplitz elements,
\[
      \tilde L_{f_k * \eta^k} {\bf T}_\phi = \tilde R_{\bar\eta^l * \bar f_l} {\bf T}_\phi = 0, 
\]
which is equivalent to the fact that the Toeplitz elements do not depend on the fibre variables $\xi,\bar \xi$.

Theorem \ref{T:covar} has the following
\begin{corollary}
The mapping  $f \mapsto {\bf T}_f$ induces an isomorphism of the algebra $(C^\infty(M)[\nu^{-1},\nu]],\star')$ onto the algebra $(C^\infty(M)[\nu^{-1},\nu]], \bullet)$ of Toep\-litz elements. 
\end{corollary}
{\it Proof.} For a function $f \in C^\infty(M)[\nu^{-1},\nu]]$ we have from Eqn. (\ref{E:Berez2}) that
\[
    {\bf T}_f = B_\star(f)\e = {\bf Q}_{B_\star(f)}.
\]
The statement of the Corollary follows from the fact that the formal Berezin transform $B_\star$ is an equivalence operator for the deformation quantizations corresponding to the star products $\star$ and $\star'$.

As it was shown in \cite{KSch}, the Berezin-Toeplitz star product on a K\"ahler manifold $(M,\omega_{-1})$ coincides with the star product $\star'$ whose opposite star product $\tilde\star$ determines the deformation quantization with separation of variables on $(M,-\omega_{-1})$ with the characterizing form $\tilde\omega$ given by Eqn. (\ref{E:dual}). Thus the construction presented in this paper can be thought of as a formal model of Berezin-Toeplitz quantization. This construction remains valid for any invertible formal constant $C(\nu)$ in the definition of the idempotent $\e$ given by Eqn. (\ref{E:defeps}).  In the rest of the section we will show that there is a natural normalization of $\e$ which determines it uniquely.

Given a formal function $f \in C^\infty(M)[\nu^{-1},\nu]]$, one can define an element ${\bf E}_f \in \E_{\nu}$ by the following local formula
\[
   {\bf E}_f = f\frac{\e}{g}\delta(\eta)\delta(\bar\eta),  
\]
where $\delta(\eta)\delta(\bar\eta)$ is the delta-function at the origin $\eta=\bar\eta=0$ so that
\[
    \int \delta(\eta)\delta(\bar\eta)d\eta d\bar \eta = 1.
\]
The fibrewise Fourier transform (\ref{E:Four}) of the element ${\bf E}_f$ is
\[
   \tilde {\bf E}_f = f\e = {\bf Q}_f.
\]
Therefore, according to Theorem \ref{T:covar}, the mapping $f \mapsto {\bf E}_f$ is a homomorphism of the algebra $(C^\infty(M)[\nu^{-1},\nu]], \star)$ to the space $\E_{\nu}$ endowed with the pullback of the product $\bullet$ via the Fourier transformation (\ref{E:Four}). If $f$ has a compact support, then, using Eqns. (\ref{E:lambda}), (\ref{E:mustar}), (\ref{E:kappa}), and (\ref{E:defeps}), we can pair the generalized function ${\bf E}_f$ with the canonical trace density $\mu_*$ of the star product $*$ as follows:
\begin{eqnarray}\label{E:pair}
    \langle {\bf E}_f,\mu_*\rangle =  \frac{\lambda_m}{\nu^{2m}}\int f \e g \, dzd\bar z = \frac{\lambda_m}{\nu^{2m}\kappa_m m!}\int f \e \omega_{-1}^m =\\
 \frac{\lambda_m}{\nu^{2m}\kappa_m m! C(\nu)}\int f\,\mu_\star.\nonumber
\end{eqnarray}
We see from Eqn. (\ref{E:pair}) that if the formal constant $C(\nu)$ is set to be
\[
           C(\nu) = \frac{\lambda_m}{\nu^{2m}\kappa_m m!},
\] 
then
\[
       \langle {\bf E}_f,\mu_*\rangle =  \int f \, \mu_\star,
\]
which means that the canonical trace density $\mu_*$ on $TM$ induces the canonical trace density $\mu_\star$ on $M$ via the mapping $f \mapsto {\bf E}_f$. Taking into account Eqns. (\ref{E:mu}) and (\ref{E:extend}) and equating the leading terms on the both sides of Eqn. (\ref{E:defeps}) we see that
\[
    \frac{1}{\nu^m m!}\omega_{-1}^m = \frac{\lambda_m}{\nu^{2m}\kappa_m m!} \nu^n \e_n \omega_{-1}^m, 
\]
whence it follows that $n = m$ and
\[
         \e_m = \frac{\kappa_m}{\lambda_m}.
\]
\section{Appendix}

Let $N$ be a submanifold of a manifold $M$ and $\I$ be the ideal of smooth functions on $M$ vanishing on $N$. We call $C^\infty(M,N) = C^\infty(M)/(\cap_r \I^r)$ the space of functions on the formal neighborhood $(M,N)$ of the submanifold $N$ in $M$. Assume that  $M$ is a complex manifold and $\bar M$ is the same manifold with the opposite complex structure. Take a local chart $U$ with holomorphic coordinates  $\{z^k,\bar z^l\}$ on $M$ and its copy $\bar U$ with coordinates 
$\{w^k,\bar w^l\}$. We will cover the diagonal $M_{\rm diag}$ of $M\times\bar M$ by the Cartesian squares $U\times\bar U$ so that on the diagonal $z^k = w^k$ and $\bar z^l = \bar w^l$. There is a mapping
\[
   \delta: C^\infty(M) \to C^\infty(M\times\bar M,M_{\rm diag})
\]
that maps a function $f(z,\bar z)$ on $M$ to its formal analytic extension $\delta f = f(z,\bar w)$ on $(M\times\bar M,M_{\rm diag})$ which is a unique solution of the equation
\[
    (\bar\p_z + \p_w)(\delta f)=0
\]
with the initial condition 
\[
\delta f|_{M_{\rm diag}} = f.
\]
Given functions $\phi,\psi\in C^\infty(M)$, we will denote by $\phi\otimes\psi$ both a function on $M\times \bar M$ and the corresponding element in  $C^\infty(M\times\bar M,M_{\rm diag})$ which will be called factorizable. Let $A$ be a differential operator on $M$. In this Appendix we will explain how to extend the bilinear operation 
\begin{equation}\label{E:operation}
    (\phi_1 \otimes \psi_1 , \phi_2 \otimes \psi_2) \mapsto (\phi_1\otimes \psi_2) \cdot \delta (A(\psi_1\phi_2))
\end{equation}
from the (linear combinations of) factorizable elements to the whole space $C^\infty(M\times\bar M,M_{\rm diag})$. A local model of $C^\infty(M\times\bar M,M_{\rm diag})$ on a chart $U\times\bar U$ can be given in the coordinates $z^k,\bar z^l, \tau^k = w^k-z^k,\bar \tau^l = \bar w^l - \bar z^l$ as 
\begin{equation}\label{E:rep}
   C^\infty(U)[[\tau,\bar \tau]],
\end{equation}
where $\tau^k,\bar\tau^l$ are treated as formal variables. Using this model one can introduce the operation
\begin{equation}\label{E:ext}
   F(z,\bar z,w,\bar w) \mapsto F|_{w=z} = F(z,\bar z, z,\bar w)
\end{equation}
on $C^\infty(M\times\bar M,M_{\rm diag})$ by setting $\tau=0$ in the formal series representing $F$ in (\ref{E:rep}).
Operation (\ref{E:ext}) is the extention of the operation
\[
   \phi\otimes\psi \mapsto (\phi\otimes 1)\cdot \delta(\psi)
\]
from the factorizable elements to the whole $C^\infty(M\times\bar M,M_{\rm diag})$. Similarly, one can introduce the operation
\[
         F \mapsto F|_{\bar z= \bar w}
\]
on $C^\infty(M\times\bar M,M_{\rm diag})$ which extends the operation
\[
   \phi\otimes\psi \mapsto \delta(\phi)\cdot(1\otimes \psi)
\]
from the factorizable elements.

Denote by $B$ the bidifferential operator on $M$ such that
\[
             B(\phi,\psi) = A(\phi\cdot\psi).
\] 
In local coordinates
\begin{eqnarray*}
    B(\phi,\psi) = B_{K\bar L P\bar Q}(z,\bar z) \left(\left(\frac{\p}{\p z}\right)^K  \left(\frac{\p}{\p\bar z}\right)^L \phi(z,\bar z)\right)\cdot\\
\left(\left(\frac{\p}{\p z}\right)^P  \left(\frac{\p}{\p\bar z}\right)^Q \psi(z,\bar z)\right).
\end{eqnarray*}
Here we assume that, say, $K = (k_1,\ldots,k_m)$ and
\[
   \left(\frac{\p}{\p z}\right)^K  = \left(\frac{\p}{\p z^1}\right)^{k_1} \ldots \left(\frac{\p}{\p z^m}\right)^{k_m},   
\]
where $m$ is the complex dimension of $M$. Now, operation (\ref{E:operation}) can be extended to the space $C^\infty(M\times\bar M,M_{\rm diag})$ as follows:
\begin{eqnarray*}
   (F_1,F_2) \mapsto B_{K\bar L P\bar Q}(z,\bar w) 
 \left(\left(\frac{\p}{\p w}\right)^K  \left(\frac{\p}{\p\bar w}\right)^L  F_1(z,\bar z,w,\bar w)\right)|_{w=z} \cdot \\
 \left(\left(\frac{\p}{\p z}\right)^P  \left(\frac{\p}{\p\bar z}\right)^Q
F_2(z,\bar z,w,\bar w)\right)|_{\bar z= \bar w}.\nonumber
\end{eqnarray*}


\begin{thebibliography}{99}
\bibitem{BFFLS}  Bayen, F., Flato, M., Fronsdal, C., Lichnerowicz, A., and Sternheimer, D.: Deformation theory and quantization. I. Deformations of symplectic structures. {\it Ann. Physics} {\bf 111} (1978), no. 1, 61 -- 110. 
\bibitem{Ber1} Berezin, F.A.: Quantization. Math. USSR-Izv. {\bf 8} (1974), 1109--1165.
\bibitem{Ber2} Berezin, F.A.: Quantization in complex symmetric spaces. Math. USSR-Izv. {\bf 9} (1975), 341--379.
\bibitem{BCG} Bertelson, M., Cahen, M., and  Gutt, S.:
Equivalence of star products. Geometry and physics.
{\it Classical Quantum Gravity}  {\bf 14}  (1997),  no. 1A, A93 -- A107.
\bibitem{BMS} Bordemann, M., Meinrenken, E., and Schlichenmaier, M.: Toeplitz quantization of K\"ahler manifolds and $gl(n),\ n \to\infty$ limits. {\it Commun. Math. Phys.} {\bf 165} (1995), 281--296. 
\bibitem{BW} Bordemann, M. and Waldmann, S.: A Fedosov star product of the Wick type for K\"ahler manifolds. {\it Lett. Math. Phys.} {\bf 41} (3) (1997), 243 -- 253.
\bibitem{CGR1} Cahen, M., Gutt S., and Rawnsley, J.: Quantization of K\"ahler manifolds I: Geometric interpretation of Berezin's quantization. {\it JGP} {\bf 7} (1990), 45--62.
\bibitem{CGR2} Cahen, M., Gutt S., and Rawnsley, J.: Quantization of K\"ahler manifolds II. {\it Trans. Amer. Math. Soc.} {\bf 337} (1993), 73--98.
\bibitem{CFS} Connes, A., Flato, M., and Sternheimer, D.: Closed star-products and cyclic cohomology, {\it Lett. Math. Phys.} {\bf 24} (1992), 1 -- 12.
\bibitem{DWL} De Wilde, M., Lecomte, P.B.A.:  Existence of
star-products and of formal deformations of the Poisson Lie algebra of
arbitrary symplectic manifolds.
{\it Lett. Math. Phys.} {\bf 7}  (1983),  no. 6, 487--496.
\bibitem{D} Deligne, P: D\'eformations de l'alg\'ebre des fonctions
d'une vari\'et\'e symplectique: comparison entre Fedosov et De Wilde, Lecomte.
{\it Selecta Math. (N.S.)} {\bf 1} (1995), no. 4, 667 -- 697.
\bibitem{F1} Fedosov, B.:  A simple geometrical construction of
deformation quantization.
{\it J. Differential Geom.} {\bf 40}  (1994),  no. 2, 213--238.
\bibitem{F2} Fedosov, B.: {\it Deformation quantization and index
theory}. Mathematical Topics, 9. Akademie Verlag, Berlin, 1996. 325 pp. 
\bibitem{Guill} Guillemin, V.: Star products on pre-quantizable symplectic manifolds. {\it Lett. Math. Phys.} {\bf 35} (1995), 85--89.
\bibitem{GR} Gutt, S. and Rawnsley, J.: Natural star products on symplectic manifolds and quantum moment maps. {\it Lett. Math. Phys.} {\bf 66}(2003), 123 --139.
\bibitem{CMP1} Karabegov, A.: Deformation quantizations with separation of variables on a K\"ahler manifold. {\it Commun. Math. Phys.} {\bf 180} (1996), 745--755.
\bibitem{Tr} Karabegov, A: On the canonical normalization of a trace density of deformation quantization,
{\it Lett. Math. Phys.} {\bf 45} (1998), 217 -- 228.
\bibitem{CMP2} Karabegov, A.: Pseudo-K\"ahler quantization on flag manifolds. {\it Commun. Math. Phys.} {\bf 200} (1999), 355--379. 
\bibitem{Deq} Karabegov, A.: On the dequantization of Fedosov's deformation quantization.
{\it Lett. Math. Phys.} {\bf 65} (2003), 133 -- 146.
\bibitem{CMP3} Karabegov, A.: Formal symplectic groupoid of a deformation quantization. {\it Commun. Math. Phys.} {\bf 258} (2005), 223--256. 
\bibitem{KSch} Karabegov, A., Schlichenmaier, M.: Identification of Berezin-Toeplitz deformation quantization. {\it J. reine angew. Math.} {\bf 540} (2001), 49-76.
\bibitem{K} Kontsevich, M.: Deformation quantization of Poisson
manifolds, I. {\it Lett. Math. Phys.} {\bf 66} (2003), 157 -- 216.
\bibitem{Mor} Moreno, C.: $*$-products on some K\"ahler manifolds. {\it Lett. Math. Phys.} {\bf 11} (1986), 361--372.
\bibitem{NT} Nest, R., Tsygan, B.: Algebraic index theorem.
{\it Commun. Math. Phys.} {\bf 172}  (1995),  no. 2, 223--262.
\bibitem{N} Neumaier, N.: Universality of Fedosov's construction for star products of Wick type on Pseudo-K\"ahler manifolds. {\it Rep. Math. Phys.} {\bf 52} (2003), 43-80. 
\bibitem{OMY} Omori, H., Maeda, Y., and Yoshioka, A.: Weyl
manifolds and deformation quantization. {\it Adv. Math.} {\bf 85} (1991), 224--255.
\bibitem{RT} Reshetikhin, N., Takhtajan, L.: Deformation quantization of K\"ahler manifolds.
L. D. Faddeev's Seminar on Mathematical Physics,  Amer. Math. Soc. Transl. Ser. 2, {\bf 201}, Amer. Math. Soc., Providence, RI, (2000), 257--276.
\bibitem{Sch} Schlichenmaier, M.: Berezin-Toeplitz quantization of compact K\"ahler manifolds. In: Quantization, Coherent States and Poisson Structures, Proc. XIV'th Workshop on Geometric Methods in Physics (Bialowieza, Poland, 9-15 July 1995), A. Strasburger, S. T. Ali, J.-P. Antoine, J.-P. Gazeau, and A. Odzijewicz, eds., Polish Scientific Publisher PWN (1998), 101 -- 115.
\bibitem{X} Xu, P.: Fedosov $*$-products and quantum momentum maps. {\it Commun. Math. Phys.} {\bf 197}  (1998),  no. 1, 167--197.
\end{thebibliography}
\end{document}